\documentclass{amsart}
\usepackage[utf8]{inputenc}
\usepackage[english]{babel}

\usepackage{forest}

\usepackage{amsmath, amssymb, amsthm, mathabx, bbm}						\usepackage{xcolor}		
\usepackage{mathtools}																						
\usepackage{placeins}
\usepackage{appendix}
\usepackage{enumerate}
\usepackage{color}
\usepackage{todonotes}

\usepackage[top=3cm, bottom=3cm, inner=2.5cm, outer=2.5cm, marginparwidth=4cm]{geometry}																		 																														 
\usepackage{booktabs}															   						
\usepackage{tikz}
\usepackage{pgfplots}
\pgfplotsset{compat=1.18}

\usepackage{enumitem}

\usepackage{float}														\usepackage{wrapfig}	

\usepackage{upgreek}										\usepackage{comment}		

\usepackage{subcaption}	
\theoremstyle{definition}

\newtheorem{theorem}{Theorem}[section]
\newtheorem{lemma}{Lemma}[section]

\newtheorem{proposition}{Proposition}[section]
\newtheorem{corollary}{Corollary}[section]
\newtheorem{claim}{Claim}[section]
\newtheorem{remark}{Remark}[section]
\numberwithin{equation}{section}
\newcommand{\floor}[1]{\left\lfloor #1 \right\rfloor}
\newcommand{\ceil}[1]{\left\lceil #1 \right\rceil}
\newcommand{\vertiii}[1]{{\left\vert\kern-0.25ex\left\vert\kern-0.25ex\left\vert #1 \right\vert\kern-0.25ex\right\vert\kern-0.25ex\right\vert}}

\begin{document}

\title[Discrete unique continuation for Schr\"odinger equations]{Discrete unique continuation for Schr\"odinger equations}

\author[T. S. Jensen, D. Pilod and J. P. G. Ramos]{Torunn S. Jensen, Didier Pilod  and Jo\~ao P. G. Ramos}

\address{Torunn S. Jensen, Department of Mathematics, University of Bergen, PO Box 7803, 5020 Bergen, Norway}
\email{Torunn.Jensen@uib.no}

\address{Didier Pilod, Department of Mathematics, University of Bergen, PO Box 7803, 5020 Bergen, Norway}
\email{Didier.Pilod@uib.no}

\address{Jo\~ao P. G. Ramos, IMPA, Rio de Janeiro, Brazil}
\email{joaopgramos95@gmail.com}

\date{\today}

\subjclass[2020]{}

\begin{abstract}
We study discrete unique continuation for one-dimensional Schrödinger equations. We prove that if a solution vanishes at the points $\pm n^{\alpha}$ at time $0$ and at the points $\pm n^{\beta}$ at time $1$, it must vanish identically in two main cases: the case of a Schrödinger equation with a real potential that is uniformly of Gevrey class $\sigma\ge 0$ in space, provided that $\max\{1,\sigma\}\big(\alpha/(1-\alpha)+\beta/(1-\beta)\big)<1$, and in the case of the cubic nonlinear Schrödinger equation, provided that $\alpha/(1-\alpha)+\beta/(1-\beta)<1$. 

Unlike previous results for logarithmically spaced zeros, our methods do not rely on complex analysis. We first show that, whenever $\alpha+\beta<1$, the zeros force the solution to lie in the Schwartz class at all times. We then convert the zeros into gaussian decay at the two times, by combining energy estimates for iterations of the vector field $\Gamma=x+2it\partial_x$ with techniques from Fourier uniqueness problems, and conclude with the Morgan-type uncertainty principle of Escauriaza, Kenig, Ponce and Vega.
\end{abstract}

\maketitle
\section{Introduction}

Most physical systems can only be observed partially: measurements are taken at a few specific times, and/or only at a discrete set of points. Schrödinger-type equations provide several examples of this:
\begin{itemize}
    \item \emph{Optical fibers.} The envelope of a light pulse obeys the cubic nonlinear Schrödinger equation, in which the distance along the fiber plays the role of time and time plays the role of the spatial variable \cite{HasegawaTappert1973, Agrawal2019}. The pulse is accessible essentially only at the two ends of the fiber, where it is recorded by detectors sampling at discrete instants.
    \item \emph{Bose--Einstein condensates.} The macroscopic wave function of Bose--Einstein condensates obeys the Gross--Pitaevskii equation, a cubic nonlinear Schrödinger equation with an external trapping potential, which reduces to a one-dimensional equation for elongated condensates \cite{DalfovoGiorginiPitaevskiiStringari1999, PitaevskiiStringari2016}. In this setting, standard absorption imaging destroys the condensate, so one obtains only a few snapshots of the object, each recorded on a pixel grid.
    \item \emph{Transmission electron microscopy.} The electron wave crossing a specimen is governed by a linear Schrödinger equation, in which the direction of propagation plays the role of time and the potential is the electrostatic potential of the specimen \cite{Kirkland2020}. The transmitted wave is recorded on a pixelated detector.
\end{itemize}
In each of these situations the detectors record the intensity $|u|^2$: in particular, this allows us to detect when the wave vanishes at a given point, if the sampled quantity also vanishes. A natural question is therefore how much of the solution is encoded in such partial information. In particular, a more direct version would be: must a solution that vanishes on a discrete set of points at two different times vanish identically?

With this question in mind, and influenced by the examples above, we study the initial value problem (IVP) associated  to the Schrödinger equation with potential
\begin{equation} \label{SP}
    \begin{cases}
        \partial_t u = i(\partial_x^2 u + V(x,t)u) & \text{in } \mathbb R \times [0,1]\\ 
        u(x,0) = u_0
    \end{cases}
\end{equation}
where $V$ is a smooth, real-valued potential satisfying some additional conditions that we specify below. We will also consider the cubic nonlinear Schrödinger equation 
\begin{equation} \label{NLS}
    \begin{cases}
        \partial_t u = i(\partial_x^2 u + \lambda|u|^2u) & \text{in } \mathbb R \times [0,1]\\ 
        u(x,0) = u_0,
    \end{cases}
\end{equation}
for some $\lambda\in \mathbb R.$
We denote $u_1=u(x,1),$ and assume that \begin{equation} \label{zero_cond} u_0(\pm n^\alpha) = u_1(\pm n^\beta) = 0 \quad \forall \, n \in \mathbb N \end{equation} for $0<\alpha,\beta<1.$
We want to understand under which conditions on $\alpha$, $\beta$ and $V$ we can deduce that $u \equiv 0$. 

More generally, assume that, for a pair of sets $(A,B)$ we have $u_0|_A=u_1|_B=0$. Our goal is to understand which conditions these sets must satisfy in order for the whole solution to be the zero solution. We refer to such a pair of sets as a \emph{Schrödinger uniqueness pair}. \\

The question above is directly motivated by the analogous case of the \emph{Fourier transform.} Indeed, in that context, we mention first the breakthrough result by Radchenko and Viazovska \cite{RadchenkoViazovska2019}, where it was shown that there are smooth functions $a_n \in \mathcal{S}_{\text{even}}(\mathbb{R})$ such that, for a given $f \in \mathcal{S}_{\text{even}}(\mathbb{R}),$ we have 
\[
f(x) = \sum_{n \ge 0} \left( f(\sqrt{n}) a_n(x) + \widehat{f}(\sqrt{n}) \widehat{a_n}(x)\right).
\]
Here, the Fourier transform is taken to have the normalization 
\[
\widehat{f}(\xi) = \int_{\mathbb{R}} f(x) e^{-2\pi i \xi x} \, dx. 
\]
The first instance in which the general question of whether a function $f \in \mathcal{S}(\mathbb{R})$ such that 
\[
f(\pm n^{\alpha}) = \widehat{f}(\pm n^{\beta}) = 0, \quad \forall n \in \mathbb{N} ,
\]
is automatically $f \equiv 0$ was posed in a paper by the third named author with M. Sousa \cite{ramos_fourier_2022}. They proved that the conclusion holds for a certain subset of the lower triangle of exponents 
\[
\left\{ (\alpha,\beta) \in (0,1)^2 \colon 0 < \alpha + \beta < 1\right\} .
\]
Subsequently, in \cite{kulikov_uniqueness}, A. Kulikov, F. Nazarov and M. Sodin showed \emph{optimal thresholds} for the question above. Effectively, they answered this question by introducing the notion of supercritical and subcritical pairs of sets. 
Indeed, they define a pair of sets $(\lambda_n)_{n\in \mathbb Z}$ and ($\mu_n)_{n\in \mathbb Z}$ to be \emph{supercritical} if $$\lim_{|n|\to +\infty}\sup |\lambda_n|^{p-1}(\lambda_{n+1}-\lambda_{n}) < \frac 12,\quad \lim_{|n|\to+\infty} \sup |\mu_n|^{q-1}(\mu_{n+1}-\mu_n) <\frac12, \quad \frac1p+\frac1q=1,$$ 
and \emph{subcritical} if 

$$\lim_{|n|\to +\infty}\inf |\lambda_n|^{p-1}(\lambda_{n+1}-\lambda_{n}) > \frac 12,\quad \lim_{|n|\to+\infty} \inf |\mu_n|^{q-1}(\mu_{n+1}-\mu_n) >\frac12, \quad \frac1p+\frac1q=1.$$ 

The main result of \cite{kulikov_uniqueness} is then that that any supercritical pair is a uniqueness pair, and that any subcritical pair is a non-uniqueness pair. This result completed the program started in \cite{ramos_fourier_2022} when it comes to showing that the Fourier uniqueness property holds in the whole lower triangle, since the pair $(\pm n^\alpha,\pm n^\beta)$ is supercritical exactly when $0<\alpha+\beta<1.$ We also refer the reader to \cite{BakanHedenmalmMontesRodriguezRadchenkoViazovska2021, Stoller2021, RamosStoller2022, CohnKumarMillerRadchenkoViazovska2022, Adve2023} for Fourier uniqueness and interpolation results in higher dimensions; to \cite{Kulikov2021, BondarenkoRadchenkoSeip2023, Lysen2025} for necessary density conditions and for uniqueness pairs in the critical regime, such as those arising from the zeros of the Riemann zeta function and other $L$-functions; to \cite{RamosSousa2023, CasseseRamos2025} for perturbations of interpolation formulas and the regularity of their basis functions; to \cite{RadchenkoStoller2022, BerghausBondarenkoRadchenkoSeipSun2025, BondarenkoSeip2026} for constructions of non-uniqueness sets; to \cite{RamosSousa2024, Lysen2026} for phaseless (Pauli-type) analogues of these questions; and to \cite{RadchenkoRamos2024} for the closely related theory of Heisenberg uniqueness pairs.

With our goal of understanding the analogous question for Schr\"odinger-like equations, we first start by noticing that a solution to the free Schrödinger equation, 
$$\begin{cases}\partial_t u = i\partial_x^2 u \\ u(x,0)=u_0, \end{cases}$$
can bee written as 
\[ u(x,t)= \frac{e^{\frac{ix^2}{4}}}{(4\pi i t)^\frac12}\left(e^{\frac{i(\cdot)^2}{4t}}u_0\right)^{\wedge} \ \Big(\frac{x}{2t}\Big), \] so that at any fixed time $t$, we may relate the solution to the Fourier transform of (a suitable modulation of) the initial data. This promptly implies that if $(A,B)$ is a \emph{Fourier uniqueness pair}, in the sense that any reasonable $f$ such that 
\[
f|_A = \widehat{f}|_B \equiv 0
\]
is automatically $f \equiv 0$, then any reasonable solution to the free Schr\"odinger equation above such that $u_0|_A = u_1|_B=0$ satisfies that $u\equiv 0.$ Since this is true for the free Schrödinger equation, it is natural to ask the question if this also holds for the solution to \eqref{SP} and \eqref{NLS}. In those cases, however, it is not immediately clear how to proceed: the relationship between the solutions to such equations and an integral transform are not always as clear as the one between the free Schr\"odinger equation and the Fourier transform, and new ideas and techniques seem to be needed. 

Indeed, progress in our main question has been restricted to \cite{kehle_uniqueness_2022}, where the authors proved that the sets $(\pm c_1 \log(1+n)^\alpha,\pm c_2 \log(1+n)^\alpha)$, for $\alpha\in(0,1/2),$ form a Schrödinger uniqueness pair for both \eqref{SP}, under suitable conditions on $V$, and \eqref{NLS}. In that setting, the techniques rely heavily on complex analysis. 

As previously mentioned, our main task is to find a suitable subset of the lower triangle of exponents where one has that $(\pm n^{\alpha}, \pm n^{\beta})$ are Schr\"odinger uniqueness pairs for either \eqref{SP} or \eqref{NLS}. Our first result is, hence, an improvement when one has certain \emph{(almost) analyticity} conditions on the potential $V$.

\begin{theorem} \label{theo:UC:potential}
    Let $V\in L^\infty([0,1]:C^\infty_b(\mathbb R))$ be real-valued and assume that there exists $A \ge 1$ and $\sigma \ge 0$ such that, for all $l \in \mathbb N$,
    \begin{equation} \label{cond:V}
    \|\partial_x^{l}V\|_{L^\infty_{x,t}}<A^{l+1} (l!)^{\sigma},
    \end{equation}
    together with
\[
\lim_{R\to+\infty}\int_0^1\sup_{|x|>R}|V(x,t)|\,dt=0.
\]
Let $u_0\in \mathcal{H}^1:=H^1(\mathbb R) \cap L^2(x^2dx)$ and $u$ be the solution to equation \eqref{SP}. Moreover assume that for $(\alpha,\beta) \in D_{\sigma} \subset \{(\alpha,\beta) \in \left(0,\frac12\right)^2 : \alpha+\beta<1\} $, $$u_0(\pm n^\alpha) =u_1(\pm n^\beta)=0 \quad \forall \, n\in\mathbb N,$$
where, in the case $0 \le \sigma \le 1$, 
\begin{equation} \label{def:Dsigma}
    D_{\sigma}= 
        \left\{ (\alpha,\beta) \in \left(0,\frac12\right)^2 : \frac{\alpha}{1-\alpha}+\frac{\beta}{1-\beta} <1 \right\} ,
\end{equation}
and, in the case $\sigma>1$, 
\begin{equation} \label{def:Dsigma.2}
    D_{\sigma}= 
        \left\{ (\alpha,\beta) \in \left(0,\frac12\right)^2 :   \frac{\alpha}{1-\alpha}+\frac{\beta}{1-\beta} <\frac1{\sigma} \right\}.
\end{equation}
Then $u\equiv 0.$
\end{theorem}

\begin{remark}
Condition \eqref{cond:V} means that $V$ is uniformly of Gevrey class $G^\sigma$ in the sense of Section~\ref{Gevrey} (see \cite{Hormander2003AnalysisI} page 281 and \cite{rodino1993linear} Section 1.4 for the classical, local Gevrey classes). More precisely, we can say the following: if $V(x,t) = V(x),$ then
\begin{itemize}
    \item[(i)] in the case $\sigma=0$, condition \eqref{cond:V} is satisfied by band-limited potentials: if $\operatorname{supp}\widehat V\subset[-\Omega,\Omega]$ and $\widehat V\in L^1(\mathbb R)$, then $\|\partial_x^lV\|_{L^\infty}\le(2\pi\Omega)^l\|\widehat V\|_{L^1}$; more generally, if $V\in L^\infty(\mathbb R)$ and $\operatorname{supp}\widehat V\subset[-\Omega,\Omega]$, then $\|\partial_x^lV\|_{L^\infty}\le(2\pi\Omega)^l\|V\|_{L^\infty}$ by the Paley--Wiener--Schwartz theorem \cite[Theorem~7.3.1]{Hormander2003AnalysisI} and Bernstein's inequality \cite{Boas1954};
    \item[(ii)] in the case $\sigma=\frac12$, condition \eqref{cond:V} is satisfied by potentials whose Fourier transform has Gaussian decay, $|\widehat V(\xi)|\le c\,e^{-\tau|\xi|^2}$ (see Lemma~\ref{lemma:G}), such as $V(x)=e^{-\tau x^2}$;
    \item[(iii)] in the case $\sigma=1$, condition \eqref{cond:V} holds for some $A>0$ if and only if $V$ extends to a bounded holomorphic function on some strip $\{z: |\operatorname{Im}z|<\delta\}$, $\delta>0$; an example is $V(x)=(1+x^2)^{-1}$;
    \item[(iv)] for any $\sigma>1$, there are smooth compactly supported functions which belong to $G^{\sigma}$ (see Proposition 1.4.10 in \cite{rodino1993linear} and Example (iv) in Subsection \ref{Gevrey} for a concrete example in our setting).
\end{itemize}
\end{remark}

\begin{figure}
\centering
\begin{tikzpicture}[scale=0.7]
\begin{axis}[
    axis lines=middle,
    xlabel={$\alpha$},
    ylabel={$\beta$},
    xmin=0, xmax=1,
    ymin=0, ymax=1,
    grid=both,
    width=10cm,
    height=10cm,
    legend style={at={(0.7,0.9)}}
]

\addplot[
    domain=0:0.5,
    samples=200,
    thick,
    blue
]
{(1 - 2*x)/(2 - 3*x)};


\addplot[
    domain=0:0.5,
    samples=200,
    fill=blue,
    fill opacity=0.2,
    draw=none
]
{(1 - 2*x)/(2 - 3*x)}
\closedcycle;
\addplot[
    domain=0:1,
    samples=2,
    thick,
    black,
    dashed
]
{1 - x};


\end{axis}
\end{tikzpicture}
\caption{The set $D_{\sigma}$, in the $\sigma = 1$ case.}
\end{figure}
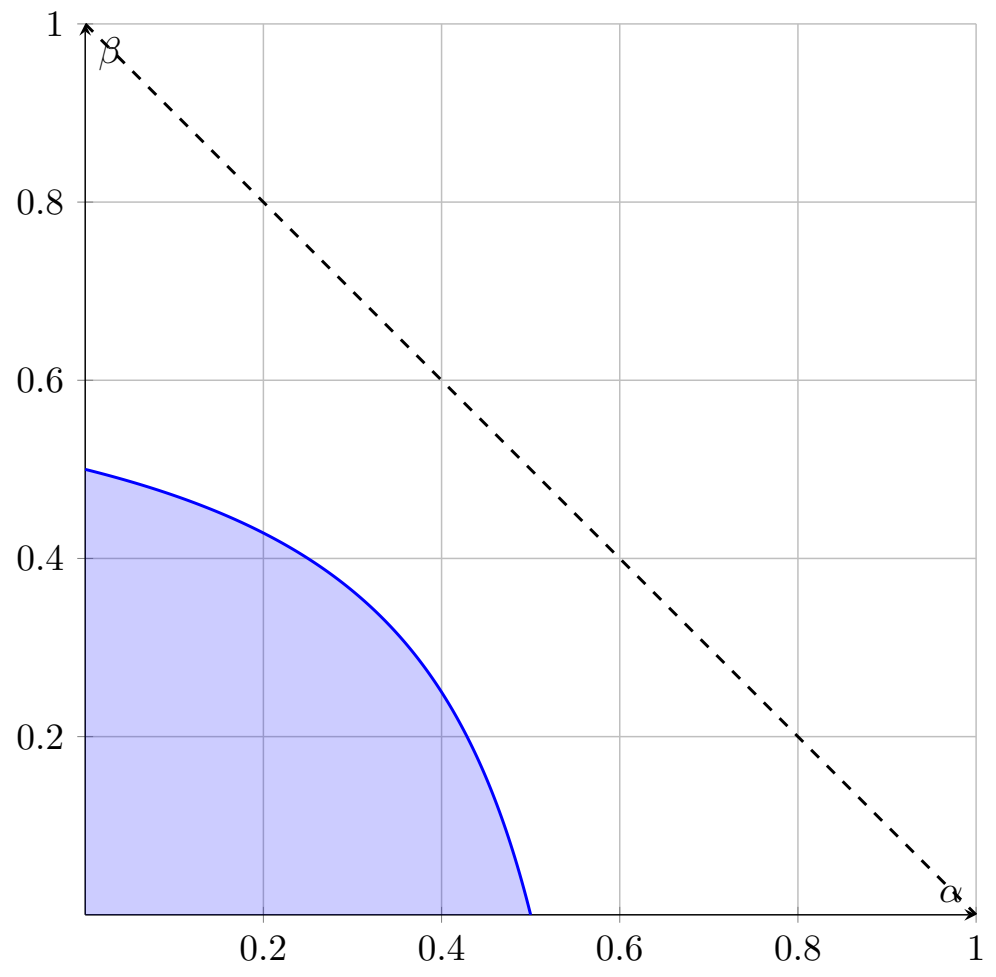

Our next result deals with the case of the \emph{cubic} nonlinear Schr\"odinger equation \eqref{NLS}: 

\begin{theorem}\label{thm:main}
Let $\lambda\in\mathbb R$ and $u_0\in \mathcal{H}^1:=H^1(\mathbb R) \cap L^2(x^2dx)$. Let $u$ be the solution to \eqref{NLS}. Assume that $u$ satisfies the zero-condition 
$$u_0(\pm n^\alpha) =u_1(\pm n^\beta)=0 \quad \forall \, n\in\mathbb N,$$
for $(\alpha,\beta)\in D_1,$ where $D_1$ is defined in \eqref{def:Dsigma.2}.
Then it follows that $u\equiv0$.
\end{theorem}

\subsection*{Strategy of the proof}
Our first step in the proof is showing, from the fact that $u_0 \in \mathcal{H}^1$ and the vanishing conditions \eqref{zero_cond}, that the solution $u$ enjoys  more regularity; more precisely, we show that $u\in C([0,1]:\mathcal{S}(\mathbb R)).$ The proof of this result builds on the results from \cite{kulikov_uniqueness} together with the fact that if a solution to \eqref{NLS} decays at two different times $u_0,u_1\in L^2(x^{2k}dx)$, then one has both regularity and decay for \emph{all} intermediate times. This latter fact relies on arguments developed in \cite{hayashi_solutions_1986,hayashi_solutions_1987,hayashi_nonlinear_1988,kehle_uniqueness_2022}.

As a matter of fact, this result may be proved in the \emph{whole} lower triangle exponent range, which is synthesized by the following results: 

\begin{theorem}\label{theo:S}
Assume that $V \in L^\infty([0,1]:C_b^{\infty}(\mathbb R))$ is real-valued.
Let $u_0\in \mathcal{H}^1:=H^1(\mathbb R) \cap L^2(x^2dx)$ and let $u$ be the solution to the equation \eqref{SP}. Moreover, assume that, for $\alpha$, $\beta \in (0,1)$ with $\alpha+\beta<1$, \begin{equation} \label{cond:u0u1}
    u_0(\pm n^\alpha)=u_1(\pm n^\beta)=0.
    \end{equation}
    Then $u\in C([0,1] : \mathcal{S}(\mathbb R))$. 
\end{theorem}

\begin{remark} Note that in this result we only require $C^{\infty}$ smoothness and boundedness of the potential: $V \in L^\infty([0,1]:C_b^{\infty}(\mathbb R))$. Moreover, the proof could be adapted in the cade of a complex-valued potential. We note furthermore, as already mentioned above, that the result actually holds in the entirety of the lower triangle. 
\end{remark}

\begin{theorem}\label{theo:S:NLS}
Let $\lambda\in\mathbb R$ and $u_0\in \mathcal{H}^1:=H^1(\mathbb R) \cap L^2(x^2dx)$ and let $u$ be the solution to  \eqref{NLS}. Moreover, assume that, for $\alpha$, $\beta \in (0,1)$ with $\alpha+\beta<1$, \eqref{cond:u0u1} holds.
    Then $u\in C([0,1] : \mathcal{S}(\mathbb R))$. 
\end{theorem}

The next step in both proofs is then showing that the zeroes actually give more information about the decay at two times. We will prove that both $$\int_{\mathbb R}e^{c|x|^p}|u_0|^2 dx \quad \text{and}\quad \int_{\mathbb R}e^{c|x|^q}|u_1|^2 dx,$$
are finite for $\frac1p +\frac1q=1$, and any $c>0.$ In order to do so, we must actually control the order of growth of norms of the form 
\[
\| x^k u_0\|_{L^2}, \, k =1, 2,\dots. 
\]
This follows essentially the same basic structure as the main theorems from \cite{ramos_fourier_2022,RamosSousa2023,kulikov_uniqueness}: one seeks to prove a bound on such quantities in terms of lower-order quantities of the same kind. In order to do so, one wants to relate this quantity to some certain \emph{dual} quantity, which in our case corresponds to bounding norms of 
\[
\Gamma_t^ju(t), \text{ where }\Gamma_t = (x+2it\partial_x). 
\]

This is where the main techniques here differ from those in \cite{kehle_uniqueness_2022}: indeed, in spite of also using the operator $\Gamma_t$, our main idea in the case of both \eqref{SP} and \eqref{NLS} is separating the \emph{main contribution} of the action of $\Gamma^k_t$ on the Duhammel expansion of $u$, from the lower-order ones. In the potential case, this is compensated by the \emph{Gevrey control} on the potential. In the cubic nonlinear case, the special structure of the nonlinearity allows for us to leave that contribution isolated, and induct on the order. 

This refined comparison estimate (see \eqref{Gamma estimate} and \eqref{est:Gammaj_nl}) is what enables us to achieve reasonable energy estimates, in contrast with what one obtains by crudely using Gr\"onwall in order to estimate the energy of $\Gamma_t^j$, as done in the proof of Theorem \ref{theo:S}. This improvement effectively means enabling the usage of the techniques from \cite{ramos_fourier_2022,kulikov_uniqueness}.

Finally, the conclusion follows by applying the Morgan principle for the Schrödinger equation with potential. We remark that the precise version below is simply one of the many statements present in the otherwise vast theory of Hardy-like uncertainty principles for Schr\"odinger equations, developed by L. Escauriaza, C. Kenig, G. Ponce and L. Vega in a series of papers \cite{Escauriaza2006Uniqueness, Escauriaza2008Convexity, Escauriaza2008Hardy, Escauriaza2010Sharp, Cowling2010Hardy, Escauriaza2011Morgan, Escauriaza2011UniqueContinuation}; see also their survey \cite{Escauriaza2012Survey}.

\begin{theorem} \label{theo:EKPV2011}\cite{Escauriaza2011Morgan}
    Given $p\in (0,1)$, there exists $N_p>0$ such that if $u\in C([0,1]:L^2(\mathbb R^n))$ is a solution of $$\partial_t u = i(\Delta u + V(x,t)u),$$ with $V=V(x,t)$ complexed valued, $\|V\|_{L^\infty(\mathbb R^n \times [0,1])}<\infty$ and 
    $$\lim_{R\to+\infty} \int_0^1\sup_{|x|>R} |V(x,t)|dt=0,$$ and there exists $\alpha,\beta>0$ such that 
$$\int_{\mathbb R^n} |u(x,0)|^2e^{2\alpha^p|x|^p/p}dx + \int_{\mathbb R^n} |u(x,1)|^2e^{2\beta^q|x|^q/q}dx<\infty,  $$
$1/p+1/q=1,$ with
$\alpha \beta >N_p,$ then $u\equiv 0.$
\end{theorem}

\begin{remark} We note that the choice of nonlinearity in Theorem \ref{thm:main} is due mostly to the fact that $2$ is the first even integer. Indeed, by following the exact same proof steps, executing the necessary modifications along the way, one is able to obtain the following: 

\begin{corollary} Let $u_0 \in \mathcal{H}^1$, and let $u$ be a solution to 
\begin{equation}\label{NLS-2k}
    \begin{cases}
        \partial_t u = i(\partial_x^2 u + \lambda|u|^{2k}u) & \text{in } \mathbb R \times [0,1]\\ 
        u(x,0) = u_0.
    \end{cases}
\end{equation}
Assume that $u$ satisfies the zero-condition \eqref{zero_cond} for $(\alpha,\beta) \in D_1.$ Then it follows that $u \equiv 0.$ 
\end{corollary}
In particular, the proof does \emph{not} rely on the integrability of the cubic nonlinear Schr\"odinger equation in dimension one in any way. In the same spirit, we also note that, by carefully redoing the same proof, one may go much further than mere powers of integers, establishing thus the following result with a gap condition, in similar terms as to that of \cite{kulikov_uniqueness}:

 \begin{corollary}\label{thm:general-zeros}
Let $(\lambda_j)_{j\in\mathbb Z}$ and $(\gamma_j)_{j\in\mathbb Z}$ be strictly increasing sequences of real numbers such that $\lambda_j\to\pm\infty$ and $\gamma_j\to\pm\infty$ as $j\to\pm\infty$. Assume that there exist $p,q>1$ such that
\begin{equation}\label{cond:gaps}
\limsup_{|j|\to\infty}|\lambda_j|^{p-1}\left(\lambda_{j+1}-\lambda_j\right)<\infty
\quad\text{and}\quad
\limsup_{|j|\to\infty}|\gamma_j|^{q-1}\left(\gamma_{j+1}-\gamma_j\right)<\infty .
\end{equation}
Let $u_0\in\mathcal{H}^1$, let $u$ be the corresponding solution of \eqref{SP} or of \eqref{NLS}. Assume that
\begin{equation}\label{cond:zeros:general}
u_0(\lambda_j)=u_1(\gamma_j)=0 \quad \text{for all } j\in\mathbb Z .
\end{equation}
\begin{itemize}
\item[(i)] Suppose that $u$ solves \eqref{SP}, where $V$ is as in the statement of Theorem \ref{thm:general-zeros}. If
\[
\frac{1}{p-1}+\frac{1}{q-1}<\frac{1}{\max\{1,\sigma\}},
\]
then $u\equiv0$.
\item[(ii)] Suppose that $u$ solves \eqref{NLS}, for any real value of the coefficient of the nonlinearity. If
\[
\frac{1}{p-1}+\frac{1}{q-1}<1,
\]
then $u\equiv0$.
\end{itemize}
\end{corollary}

We have, however, chosen to keep the power of integer case, as well as the cubic NLS case, for simplicity of exposition. 
\end{remark}

The paper is organized as follows. In Section \ref{sec:prelim}, we recall Stirling's formula, describe the uniform Gevrey class $G^\sigma$, and recall some identities related to the vector field $\Gamma^k_t$. In Section \ref{sec:S}, we explain how an initial condition in $\mathcal{H}^1$ together with the zero condition at two different times imply that the solution lies in the Schwartz class for all time; in short, we give the proof of Theorems \ref{theo:S} and \ref{theo:S:NLS}. Finally, Section \ref{sec:UC_potential} is dedicated to the proof of the unique continuation theorem in the case of the potential (Theorem \ref{theo:UC:potential}), while Section \ref{sec:UC_NLS} is devoted to the proof of the unique continuation theorem in the case of the nonlinear Schr\"odinger equation.

\subsection*{Notation} Throughout the text, $c, c_\beta, c_\alpha$ will denote arbitrary positive constants, only depending on $\alpha$ and $\beta$, and that may change from line to line. In Section \ref{sec:S}, the constants $c$ may also depend on the indexes $k, j, l$ involved in the proof, since the growth of these do not matter at this point. However, in Sections \ref{sec:UC_potential} and \ref{sec:UC_NLS}, the constants $c$ only depend on $\alpha$ and $\beta$, and we keep track carefully on the growth of the indexes $k, j, l$.

\subsection*{Acknowledgements} J. P. G. R.  was supported by the FCT through project SHADE
(2023.17881.ICDT), DOI \allowbreak \texttt{10.54499/}
\texttt{2023.17881.ICDT}. He was also supported by FAPERJ through  JCNE grant
no.~SEI-260003/ \allowbreak 020475/2025, by Instituto Serrapilheira through
grant Serra-R-2510-59765, and by CNPq through Research Productivity Fellowship no. 303861/2026-6. T.S.J. and D.P. were partially supported by the Trond Mohn Foundation. They also would like to thank the University of Lisbon and the Insituto de Matm\'atica Pura e Applicada (IMPA) for the kind hospitality during different stages of this work.

\subsection*{LLM Usage.} All the ideas and arguments, as well as their technical details, were conceived and verified by the authors without using artificial intelligence tools. In particular, in Sections \ref{sec:S} and \ref{sec:UC_potential}, LLM were only used to check computations and language. A preliminary working version of Section \ref{sec:UC_NLS} was typeset by using several Large Language Models, but the present final version was completely written by the authors. Large Language Models were only used in the final stages for looking for typos and searching references.

\section{Preliminaries} \label{sec:prelim}
\subsection{Stirling formula}
We recall the definition of the Gamma function\footnote{We also use $\Gamma_t$ to denote the vector field $x+2it\partial_x$ (see Subsection \ref{subsec:gamma}). It will however be clear to the reader which one that should be considered.}
$$\Gamma(t) = \int_0^{+\infty} e^{-x} x^{t-1}dx, \quad t>0.$$
For integers $k\ge 0,$
$$\Gamma(k+1)=k!.$$ Moreover, the Stirling formula states that
\begin{equation} \label{Stirling_Gamma} \Gamma(t+1) = t^{t}e^{-t}\sqrt{2\pi t} \left(1+O\Big(\frac1t\Big)\right), \quad \text{as} \ t \to +\infty.\end{equation}
In particular, \begin{equation} \label{Stirling_Factorial} k!=k^ke^{-k}\sqrt{2\pi k}\left(1+O\Big(\frac1k\Big)\right), \quad \text{as} \ k \to +\infty.\end{equation} 
As a consequence, we will often use the simpler identities $k!=e^{k\log k + O(k)}$ and $\Gamma(t) = e^{t\log{t}+ O(t)}$, which provide the useful bounds
    \begin{equation}\label{bound:factorial}
   k!\le  k^k \le e^k k!.
    \end{equation}

\subsection{Gevrey Spaces} \label{Gevrey}
We recall the definition of Gevrey functions of order $\sigma \ge 0$, (see \cite{Hormander2003AnalysisI} and Chapter I in \cite{rodino1993linear}),
$$  \left\{ V \in C^\infty(\mathbb{R}) : |\partial_x^j f(x)| \le c_K^{j+1}(j!)^\sigma, \quad\forall \, x\in K,\ j\in\mathbb{N} \right\}, $$
for any compact $K\subset \mathbb R.$ and $c_K$ is a positive constant depending on $K.$

In this paper, we work with a stronger assumption on the potential.  We say that $V$ is \emph{uniformly of Gevrey class} $\sigma \ge 0$ if there exists $A>0$ such that 
\[\sup_{x\in \mathbb R} |\partial_x^j V(x)|\le A^{j+1} (j!)^\sigma,  \forall \, j\in \mathbb N.\] 
If $V$ is uniformly of Gevrey class $\sigma$, we denote $V \in G^\sigma=G^\sigma(\mathbb R)$. Note that the notation $G^\sigma$ differs from \cite{Hormander2003AnalysisI,rodino1993linear}, since it means uniformly of Gevrey class on $\mathbb R$ in our case. 

In the following, we derive some useful properties of the spaces $G^\sigma(\mathbb R)$ and give some important examples. The rest of this subsection is largely inspired by Chapter 1 in \cite{rodino1993linear}.  

Observe that $G^{\sigma_1} \subset G^{\sigma_2}$ if $\sigma_1\le \sigma_2$. 

We start with a characterization of the spaces $G^\sigma$ by using the Fourier transform.
\begin{lemma} \label{lemma:G}
Suppose that there exist $s >0$, $\tau>0$ and $c>0$ such that $|\hat{f}(\xi)| \le ce^{-\tau |\xi|^{s}}$,  then $f \in G^\sigma$, for $\sigma=1/s$. In particular, the Gaussian $e^{-\tau x^2}$ is uniformly of Gevrey class of order $\sigma=1/2.$
\end{lemma}

\begin{proof}
Observe that
\begin{align*}
    \|\partial_x^j f\|_{L^\infty}&\le \int_{\mathbb R} |\xi|^j \hat{f} 
(\xi) d\xi \le \int_{\mathbb R} |\xi|^j ce^{-\tau |\xi|^s} d\xi .\end{align*}
We consider the integral $$I_{j,s}=\int_{0}^{+\infty} r^j e^{-\tau r^s} dr.$$
By doing the change of variable $r=\left(\frac{t}{\alpha}\right)^{1/\beta},$ it follows that

\begin{align*}
    I_{j,s}= s^{-1} \tau^{-(j+1)/\beta}\int_{0}^{+\infty} t^{(j+1)/s-1}e^{-t} dt = s^{-1} \tau^{-(j+1)/s} \Gamma\left(\frac{j+1}{s}\right) \le c_{s,\tau}^j e^{\frac js\log {j}+ O(j)} \le c_{s,\tau}^j (j!)^{1/s}.
\end{align*}
by using the Stirling formula for the $\Gamma$ function in \eqref{Stirling_Gamma}. 

Thus, if $|\hat{f}(\xi)| \le c e^{-\tau |\xi|^{s}}$, then $\|\partial_x^j f\|_{L^\infty} \le c_{s,\tau}^j (j!)^{1/s}$, which means that $f\in G{^{1/s}}.$ 
\end{proof}

\noindent {\bf Examples of uniformly Gevrey potentials depending on $\sigma \ge 0$}.

\medskip
\noindent (i) \emph{$\sigma=0$}. If $V\in L^\infty(\mathbb R)$ is band-limited, \textit{i.e.} $\text{supp} \, \widehat{V} \subset [-A,A]$, for some $A>0$, then $V \in G^0$.

\medskip
\noindent (ii) \emph{$1/2 \le \sigma <1$}. Let $m\ge1$ be an integer, let $\tau>0$, and let $V(x)=e^{-\tau x^{2m}}$. Then $V$ extends to an entire function which satisfies
\[
|V(x+iy)|\le e^{-\frac{\tau}{2}x^{2m}+C_m\tau|y|^{2m}}, \qquad x,y\in\mathbb R,
\]
for some constant $C_m>0$ depending only on $m$. By the characterization of the Gelfand--Shilov spaces in terms of entire extensions \cite[Chapter IV]{GelfandShilov1968vol2}, this means that $V$ belongs to the Gelfand--Shilov space $S^{1-\frac1{2m}}_{\frac1{2m}}$. Since the Fourier transform maps $S^{\beta}_{\alpha}$ onto $S^{\alpha}_{\beta}$ (see \cite[Chapter IV]{GelfandShilov1968vol2} and \cite{ChungChungKim1996}), it follows that
\[
|\widehat{V}(\xi)|\le c\,e^{-\tilde{\tau}|\xi|^{\frac{2m}{2m-1}}}
\]
for some $c, \tilde{\tau}>0$ (see also \cite{Boyd2014quartic} for the precise asymptotic behaviour of $\widehat V$). Lemma \ref{lemma:G} then implies that $V\in G^{\frac{2m-1}{2m}}$. In particular, for any $\tau>0$, $e^{-\tau x^2} \in G^{\frac12}$ and $e^{-\tau x^{4}}\in G^{\frac34}$.

\medskip 
\noindent (iii) \emph{$\sigma =1$}. The class $G^1$ corresponds to analytic functions on $\mathbb R$ that can be extended holomorphically at least on the strip $\{z:  |\text{Im} (z)| <1/A\}$.  A good example is the Poisson kernel $P(x)=\frac1{1+x^2}$ whose Fourier transform is $\widehat{P}(\xi)=e^{-2\pi |\xi|}$. So, it follows from Lemma \ref{lemma:G} that $P \in G^1$.

\medskip
\noindent (iv) \emph{$\sigma>1$.} We will show that there exist compactly supported potentials belonging to $G^\sigma$, $\sigma>1$.

For any $s>0$, 
let us define  $$V(x)=\begin{cases}
    e^{-\frac{1}{(1-x^2)^s}} & |x|<1, \\ 0 & |x|\ge 1 .
\end{cases}$$
Then, it is well-known that $V \in C_c^\infty(\mathbb R)$. Moreover, we claim that $\|\partial_x^j V\|_{L^\infty} \le c^j (j!)^{(1+\frac1s)}$, so that $V\in G^{1+\frac1s}.$ Of course this does not work for $G^1$, since there are no nonzero analytic compactly supported functions. To show the claim, we compute the derivatives: $V'(x) =-e^{-\frac1{(1-x^2)^s}} \frac{2xs}{(1-x^2)^{s+1}}$.
Since $|x|<1$ in $\text{supp} \, V,$  and by denoting  $t=1-x^2 \in [0,1]$, it follows that
$$|V'(x)|\le e^{-\frac1{(1-x^2)^s}} \frac{2s}{(1-x^2)^{s+1}} \le c e^{-1/t^s}t^{-(s+1)}.$$
Similarly, we can show that
$|V^{(j)}(x)|\le c^je^{-1/t^s} t^{-(s+1)j}\le c^j e^{-\frac1{t^s} -(s+1)j\log t}$ for any $j\ge 1$.
Moreover,
$$\frac d{dt}(-\frac1{t^s} -(s+1)j \log t) = \frac s{t^{s+1}}-(s+1)j \frac 1t =0 \iff t^s=\frac s{(s+1)j},$$
and the function $f(t)=-\frac1{t^s} - (s+1)j \log t$ attains its maximum for $t=\left(\frac s{(s+1)j}\right)^{1/s}$.
It follows that
$$|V^{(j)}(x)|\le c^je^{-\frac{(s+1)j}{s} -(s+ 1)j \log \left(\frac s{(s+1)j}\right)^{1/s}}\le c^j e^{ \frac{(s+1)}{s}j \log j } \le c^j (j!)^{ 1+\frac1s},$$
which shows that $V\in G^{1+\frac1s}$.

\subsection{The operator $\Gamma_t$} \label{subsec:gamma}
We recall from \cite{ hayashi_solutions_1986,hayashi_solutions_1987,hayashi_nonlinear_1988} the operator 
\begin{equation} \label{Gamma:form.1}
    \Gamma_t^j = (x+2it\partial_x)^j .
\end{equation}
It is known that for $j\in \mathbb{N}$ and $t\ne 0,$
\begin{equation} \label{Gamma:form.2}
 \Gamma^j_t u(t) = e^{i x^2/4t}(2it)^{j}\partial_x^j(e^{-i x^2/4t}u(t)),   
\end{equation}
and that \begin{equation}\label{Gamma:eit}
     \Gamma^j_t f = e^{it\partial_x^2}\left(x^j e^{-it\partial_x^2} f\right).
\end{equation}
We let, for $t\ne 0,$ $\tilde{u}(x,t)=e^{-i|x|^2/4t}u(x,t),$ so that 
\begin{equation} \label{Gamma:form.3}
\Gamma_t^j u= e^{i|x|^2/4t} (2it)^{|j|}\partial_x^{j} \tilde{u} \quad \text{and} \quad \|\Gamma^j_t u\|_{L^2} = (2t)^{|j|}\|\partial_x^j \tilde u\|_{L^2}.
\end{equation}

Moreover $\Gamma^j_t$ commutes with the free Schrödinger flow, that is
\begin{equation}
    \Gamma^j_t (i\partial_t + \partial_x^2) = (i\partial_t + \partial_x^2) \Gamma^j_t,
\end{equation}
 and if $u$ is a solution to \eqref{SP}, then $\Gamma^j_t u $ is a solution to the IVP
\begin{align}\label{equation for Gamma u}
    \begin{cases}
        \partial_t \Gamma^j_t u = i(\partial_x^2\Gamma^j_t u + \Gamma^j_t (Vu)) \\
        \Gamma^j_0u(x,0) = x^j u_0.
    \end{cases}
\end{align}

Since the Schrödinger flow is reversible in time,
the solution at time $t=0, \ u_0,$ can be recovered from the data at time $t=1$.
We can therefore write $u_0$ as the solution obtained by evolving $u_1$
backwards through the same equation. More precisely, if we let $v(x,t)=u(x,1-t),$ then 
$$\begin{cases}\partial_t v = -i(\partial_x^2 v + V(x,1-t)v) \\ v(x,0)=u_1,
\end{cases}$$
and $u_0 = v(x,1)=v_1.$
Then we can work with the operators
$$\Gamma_{-t}^j =(x-2it\partial_x)^j ,$$
so that
\begin{equation} \label{Gamm:form.-1}
     \Gamma^j_{-t} v(t) = e^{-i x^2/4t}(-2it)^{j}\partial_x^j(e^{i x^2/4t}v(t))
\end{equation}
and
\begin{equation} \label{Gamm:form.-2}\Gamma^j_{-t} f= e^{-it\partial_x^2}\left(x^j e^{it\partial_x^2}f\right). 
\end{equation}
We define for $t\ne 0$ $\tilde{v}_-:=e^{i x^2/4t}v(t),$
so that $\tilde{v}_{-,1} = e^{ix^2/4} v_1 = e^{ix^2/4}u_0.$ Then 

$$\Gamma_{-t}^j v=e^{-i x^2/4t}(-2it)^{j} \partial_x^j \tilde{v}_-, \quad \text{ and } \quad \|{\Gamma}^j_{-t} v \|_{L^2}=(2t)^{j}\|\partial_x^j \tilde{v}_-\|_{L^2}.$$
Moreover, ${\Gamma_{-t}^j}$ commutes with the \lq\lq backwards\rq\rq \,  Schrödinger flow, so that $\Gamma^j_{-t} v(t)$ satisfies the IVP
\begin{align}\label{equation for Gamma v}
    \begin{cases}
        \partial_t \ \Gamma^j_{-t} v = -i(\partial_x^2\Gamma^j_{-t} u + \Gamma^j_{-t} (V(\cdot,1-t)v)) \\
         \Gamma^j_0v(x,0) = x^j v_0 = x^j u_1.
    \end{cases}
\end{align}

\section{From $\mathcal{H}^1$ to $ \mathcal{S}(\mathbb{R})$.} \label{sec:S}

For $k\in \mathbb{N}$ we define
\[ \mathcal{H}^k:=H^k(\mathbb R) \cap L^2(\mathbb R;x^{2k}dx),\] which is a Banach space with the norm 
$\|u\|_{\mathcal{H}^k}^2 = \sum_{l=0}^{k} \left( \|\partial_x^l u\|_{L^2}^2 +\|x^l u\|_{L^2}^2\right)$. 
We recall the following estimate whose proof can be found in Lemma 4.3 in \cite{kehle_uniqueness_2022}, 
\begin{equation} \label{interp:Hk}
\sum_{l+j\le k} \|x^l \partial_x^j f\|_{L^2} \lesssim \sum_{l,j \le k} \left(\|x^l f\|_{L^2} + \|\partial_x^j f\|_{L^2} \right).
\end{equation}
The above estimate shows that one may use the equivalent norm $\sum_{l+j\le k}\|x^{l}\partial_x^j u\|_{L^2}$ on $\mathcal{H}^k$. 

We also define
\[C_b^{k}(\mathbb R):= \left\{V \in C^{k}(\mathbb R : \mathbb R) : \sup_{x \in \mathbb R}|V^{(j)}(x)|<\infty, \, \forall \, j\le k\right\}\]
and
\[C_b^{\infty}(\mathbb R):= \left\{f \in C^{\infty}(\mathbb R : \mathbb R) : \sup_{x \in \mathbb R}|V^{(j)}(x)|<\infty, \, \forall \, j \in \mathbb N\right\}.\]

We start by stating a well-posedness result for the initial value problem \eqref{SP} in $\mathcal{H}^k(\mathbb R)$ when $V\in L^\infty([0,1]:C^k_b(\mathbb R)).$
\begin{proposition} \label{WP}
    Assume that $V\in L^\infty([0,1]:C^k_b(\mathbb R))$ and let $k \ge 1$. Then, for any initial datum $u_0 \in \mathcal{H}^k$, there exists a unique solution to \eqref{SP} in $C([0,1] : \mathcal{H}^k)$. Moreover, the flow map $:u_0 \in \mathcal{H}^k \mapsto u \in C([0,1] : \mathcal{H}^k)$ is Lipschitz.
\end{proposition}

\begin{proof} The proof is a classical application of the fixed point theorem, and follows from the arguments in the proof Theorem 1 in \cite{hayashi_nonlinear_1988}.

\end{proof}

Since $\displaystyle{\mathcal S(\mathbb R)=\bigcap_{k=1}^{+\infty} \mathcal{H}^k(\mathbb R)},$ we get the following corollary.

\begin{corollary} \label{W.p.S}
  Assume that $V\in L^\infty([0,1]:C^\infty_b(\mathbb R))$. Then, for any initial data $u_0 \in \mathcal S(\mathbb R)$, there exists a unique solution to \eqref{SP} in $C([0,1] : \mathcal S(\mathbb{R}))$. Moreover, the flow map $:u_0 \in \mathcal{S}(\mathbb R) \mapsto u \in C([0,1] : \mathcal{S}(\mathbb R ))$ is Lipschitz.
\end{corollary}

Now, we state a similar well-posedness result in $\mathcal{H}^k$ for the nonlinear Schr\"odinger equation \eqref{NLS}.

\begin{proposition} \label{prop:LWP:NLS}
Let $k \ge 1$. For any initial datum $u_0\in\mathcal{H}^k$, there exists a unique solution $u$ to \eqref{NLS} in $C([0,1]:\mathcal{H}^k)$ satisfying 
\begin{equation} \label{LWP:NLS:bound}
\sup_{ 0 \le t \le 1} \|u(t)\|_{\mathcal{H}^k}\le c \|u_0\|_{\mathcal{H}^k} .
\end{equation}
Moreover, the flow map $:u_0 \in \mathcal{H}^k \mapsto u \in C([0,1] : \mathcal{H}^k)$ is Lipschitz.
\end{proposition}

\begin{proof} The proof is a classical application of the fixed point theorem on the Duhamel formulation associated to \eqref{NLS}, following the arguments in the proof of Theorem 1 in \cite{hayashi_nonlinear_1988}. In order to extend the existence time to $[0,1]$, we use tamed estimates on the Duhamel formulation of \eqref{NLS} to obtain \text{a priori} estimates on the $\|u(t)\|_{\mathcal{H}^k}$-norm on a time existence depending only on the $\|u(t)\|_{\mathcal{H}^1}$ norm, and then we use the $\mathcal{H}^1$ control on the solutions obtained in Theorem 2.1 in \cite{hayashi_solutions_1986} to conclude. 
\end{proof}

As a corollary, we state the well-known well-posedness result for \eqref{NLS} in $\mathcal{S}(\mathbb R))$. 

\begin{corollary} \label{W.p.S.NLS}
  For any initial data $u_0 \in \mathcal S(\mathbb R)$, there exists a unique solution to \eqref{NLS} in $C([0,1] : \mathcal S(\mathbb{R}))$. Moreover, the flow map $:u_0 \in \mathcal{S}(\mathbb R) \mapsto u \in C([0,1] : \mathcal{S}(\mathbb R ))$ is Lipschitz.
\end{corollary}

As a consequence of Proposition \ref{WP}, we prove that the decay of a solution at two different times imply regularity and decay for all times in between.
\begin{proposition} \label{Prop2} Assume that $V\in L^\infty\left([0,1] : C^k_b(\mathbb{R})\right).$
    Let $u_0 \in \mathcal{H}^1$ and let $u$ be the solution $u\in C([0,1]: \mathcal{H}^1)$ obtained in Proposition \ref{WP}. Moreover, assume that $u_0,u_1=u(\cdot,1)\in L^2(x^{2k}dx) $ for some $k\in \mathbb N$. Then $u\in C([0,1]:\mathcal{H}^{k}(\mathbb R)).$
\end{proposition}
\begin{proof}

Let $\{u_0^m\}_m$ be a sequence of Schwartz functions such that $u_0^m\to u_0$ in $\mathcal{H}^1 \cap L^2(x^{2k},dx)$. 
Let $\{u^m\}\subset C[0,1] : \mathcal{S}(\mathbb R))$ be the sequence of the corresponding solutions obtained in Corollary \ref{W.p.S}. By Proposition \ref{WP} it follows that $u^m\to u$ in $C([0,1]:\mathcal{H}^1) \subset C([0,1]:L^2(\mathbb R))$. 

Now $\Gamma^k_t u^m$ is well defined, but we need to make sense of $\Gamma^k_t u$. The first step is to prove that $\Gamma^k_t$ can be defined as $$\Gamma^k_t u:=\lim_{m\to +\infty} \Gamma^k_t u^m$$ and that $\Gamma^k_t u\in C([0,1]:L^2(\mathbb R))$. 

By making an energy estimate, we have 
\begin{align}
 \frac{1}{2} \frac{d}{dt} \|\Gamma^k_t u^m(t)\|_{L^2}^2 &= Re \int_{\mathbb{R}}\partial_t \Gamma^k_t u^m(t) \overline{\Gamma^k_tu^m(t)} dx \nonumber\\ & = Re \left( i\int_{\mathbb{R}} \partial_x^2 \Gamma^k_t u^m(t) \overline{\Gamma^k_tu^m(t)} dx \right) + Re \left( i \int_{\mathbb{R}} \Gamma^k_t(Vu^m(t)) \overline{\Gamma^k_t u^m(t)} dx\right)\nonumber \\ & = -Im\left( \int_{\mathbb{R}} \Gamma^k_t(Vu^m(t)) \overline{\Gamma^k_t u^m(t)} dx\right) \nonumber \\ & = -Im\left( (2t)^{2k} \int_{\mathbb{R}}\partial_x^k(V\tilde{u}^m(t))\overline{\partial^k_x \tilde{u}(t)}dx \right) \nonumber\\ & = -\sum_{l=0}^{k-1}\binom{k}{l}Im\left((2t)^{2k} \int_{\mathbb{R}} \partial_x^{k-l}V \partial_x^l \tilde{u}^m(t) \overline{\partial_x^k \tilde{u}^m(t)} dx\right) + Im \int_{\mathbb{R}}V |\partial_x^k \tilde{u}^m(t)|^2dx\nonumber\\ & = - \sum_{l=0}^{k-1}\binom{k}{l}Im\left((2t)^{2k} \int_{\mathbb{R}} \partial_x^{k-l}V \partial_x^l \tilde{u}^m(t) \overline{\partial_x^k \tilde{u}^m(t)} dx\right) \label{ener:est}\\ &\le c(k)\sup_{0\le l\le k}\|\partial_x^{l} V\|_{L^\infty}\sum_{l=0}^{k-1}\|\partial_x^l \tilde{u}^m\|_{L^2}\|\partial_x^k \tilde{u}^m\|_{L^2}, \nonumber
\end{align}
so that, by using our assumptions on $V$, 
\begin{equation*}
    \frac d{dt}\|\Gamma^k_t u^m(t) \|_{L^2} \le c(k,V) \sum_{l=0}^{k-1}\|\partial_x^l \tilde{u}^m(t)\|_{L^2} \le c(k,V)(\|\tilde{u}^m\|_{L^2} + \|\partial_x^k \tilde{u}^m(t)\|_{L^2}) \le c(k,V)(\|u_0^m\|_{L^2} + \|\Gamma^k_tu^m(t)\|_{L^2}). 
\end{equation*}
Here we used that for each $0\le l\le k-1$ we can write
$$\|\partial_x^l \tilde{u}^m\|_{L^2}^2 = \int_{|\xi|\le 1}\xi^{2l}|\tilde{u}^m|^2 d\xi +\int_{|\xi|\ge1} \xi^{2l}|\tilde{u}^m|^2 \le \|u_0^m\|_{L^2}^2 + \|\partial_x^k \tilde{u}^m\|_{L^2}^2.  $$
Thus, by Grönwall's inequality, we get that, for all $0 \le t \le 1$,
\begin{equation}
    \|\Gamma^k_t u^m(t)\|_{L^2} \le e^{c(k,V)t} \left(\|\Gamma^k_0 u_0^m\|_{L^2} + \|u_0^m\|_{L^2}\right) \le c(k,V)\left(\|x^k u_0^m\|_{L^2} + \|u_0^m\|_{L^2}\right),
\end{equation}
where the right hand side is uniformly bounded. Thus
\begin{equation} \label{sup_m Gamma_k,m}
    \sup_m \sup_{t\in [0,1]}\|\Gamma^k_t u^m(t)\|_{L^2} \le c(k,V)\sup_m \left(\|x^k u_0^m\|_{L^2} + \|u_0^m\|_{L^2}\right).
\end{equation}
Since $u_0^m\to u_0$ in $\mathcal{H}^1\cap L^2(x^{2k}dx),$ we know that the right hand side of \eqref{sup_m Gamma_k,m} converges to 
$c\left(\|x^k u_0\|_{L^2} + \|u_0\|_{L^2}\right).$ 

By the same estimate as above, we have for all $t\in[0,1]$,
\begin{align*}
 \|\Gamma^k_t u^{m_1}(t)-\Gamma^k_t u^{m_2}(t)\|_{L^2}  \le c(k,V)\left(\|x^k (u_0^{m_1}-u_0^{m_2})\|_{L^2} + \|u_0^{m_1}-u_0^{m_2}\|_{L^2}\right)\underset{m_1,m_2\to +\infty}{\longrightarrow} 0,
\end{align*}
which shows that $\{\Gamma^k_t u^m\}_m$ is a Cauchy sequence in $C([0,1]:L^2(\mathbb{R})).$ Since $C([0,1] : L^2(\mathbb{R}))$ is complete, there exists a unique $w$ such that $\Gamma^k_t u^m\to w$ in $C([0,1] : L^2(\mathbb{R})).$  

Moreover, for each $m$ we have that 
$\Gamma^k_t u^m = e^{ix^2/4t}(2it)^{k}\partial_x^k \tilde {u}^m.$ Then, for all $t \in (0,1]$,
$$\partial_x^k \tilde{u}^m \underset{n \to +\infty}{\longrightarrow} h=e^{-ix^2/4t}(2it)^{-k}w,  \quad \text{in} \ L^2(\mathbb R)$$ and in particular $h=\partial_x^k \tilde{u}$ in the weak sense, because $\{u^m\}$ converges to $u$ in $C([0,1] : L^2(\mathbb R))$. Moreover
$$ \Gamma^k_t u^m =(2it)^k e^{ix^2/4t} \partial_x^k \tilde{u}^m \rightharpoonup (2it)^ke^{ix^2/4t}\partial^k_x \tilde{u} \quad \text{weakly in} \ L^2(\mathbb R).$$Since $\Gamma^k_t u^m \to w$ strongly, we also have
$\Gamma^k_t u^m \rightharpoonup w$, and by uniqueness of the weak limit, we must have $\Gamma^k_t u:=w = (2it)^ke^{ix^2/4t}\partial^k_x \tilde{u}.$
In particular we get the bound
\begin{equation} \label{bound on Gamma}
    \sup_{t\in [0,1]}\|\Gamma_t^k u\|_{L^2} \le c(k,V)\left(\|x^k u_0\|_{L^2} + \|u_0\|_{L^2}\right).
\end{equation}
The next step is to show that $u_1$ belongs to $\mathcal{H}^k$. Since $\|\Gamma^j_1 u(1)\|_{L^2}=2^j\|\partial_x^j \tilde{u}_1\|_{L^2}$ and $\|x^j\tilde{u}_1\|_{L^2}=\|x^j u_1\|_{L^2},$ it follows that $\partial_x^j \tilde{u}_1\in L^2$ and $x^j \tilde{u}_1\in L^2$ for all $j\le k.$ This implies that $\tilde{u}_1 \in \mathcal{H}^k$.

By using estimate \eqref{interp:Hk}, we deduce that 
$$\sum_{l+j\le k} \|x^l\Gamma_1^j u_1\|_{L^2} = \sum_{l+j\le k} 2^j\|x^l \partial_x^j \tilde{u}_1\|_{L^2}\lesssim \|\tilde{u}_1\|_{\mathcal{H}^k} .$$ 

Finally, since $\Gamma^j_1=(x+2i\partial_x)^j,$ it follows by induction that
$$\partial_x^j = \sum_{l+s \le j}c(l,s,j)x^l \Gamma^s_1,$$
for some positive constants $c(l,s,j)$. This implies that
$$\sum_{l+j\le k}\|x^l \partial_x^j u_1\|_{L^2}\lesssim \|\tilde{u}_1\|_{\mathcal{H}^k},$$
so that $u_1\in \mathcal{H}^k.$

We can prove a similar estimate to \eqref{bound on Gamma} for the function $v(x,t)=u(x,1-t)$, 
\begin{equation}\label{Bound on Gamma reverse}
    \sup_{t\in [0,1]}
\|\Gamma^k_{-t}v(t)\|_{L^2} = \sup_{t\in [0,1]}\|\Gamma^k_{t-1}u(t)\|_{L^2} \le c(k,V)\left(\|x^k u_1\|_{L^2} + \|u_0\|_{L^2}\right).\end{equation}
which by arguing as for $u_1$ will imply that $u_0\in \mathcal{H}^k.$ By Proposition \ref{WP}, we finally conclude that $u\in C([0,1] : \mathcal{H}^k).$
\end{proof}

\begin{corollary}\label{Cor. S}
Assume that $V\in L^\infty([0,1]:  C^\infty_b(\mathbb R))$.   Let $u_0 \in \mathcal{H}^1$ and let $u$ be the solution $u\in C([0,1]: \mathcal{H}^1)$ obtained in Proposition \ref{WP}. Moreover, assume that $u_0,u_1=u(\cdot,1)\in L^2(x^{2k}dx) $ for all $k\in \mathbb N$. Then $u\in C([0,1]:\mathcal S(\mathbb R))$.
\end{corollary}

 Similar results are known for solutions to \eqref{NLS}. The proof is almost identical to the proof of Proposition \ref{Prop2}, see \cite{hayashi_solutions_1986, hayashi_solutions_1987, hayashi_nonlinear_1988, kehle_uniqueness_2022}.
 \begin{proposition}
    Let $u_0 \in \mathcal{H}^1$ and let $u$ be the solution $u\in C([0,1]: \mathcal{H}^1)$ obtained in Proposition \ref{prop:LWP:NLS}. Moreover, assume that $u_0,u_1=u(\cdot,1)\in L^2(x^{2k}dx) $ for some $k\in \mathbb N$. Then $u\in =C([0,1]:\mathcal{H}^{k}).$
 \end{proposition}
 \begin{corollary}\label{Cor. S:NLS}
Let $u_0 \in \mathcal{H}^1$ and let $u$ be the solution $u\in C([0,1]: \mathcal{H}^1)$ obtained in Proposition \ref{prop:LWP:NLS}. Moreover, assume that $u_0,u_1=u(\cdot,1)\in L^2(x^{2k}dx) $ for all $k\in \mathbb N$. Then $u\in C([0,1]:\mathcal S(\mathbb R))$.
\end{corollary}

Now, we are in position to prove Theorem \ref{theo:S} and Theorem \ref{theo:S:NLS}.
\begin{remark} 
Thanks to Proposition \ref{WP}, the corresponding solution $u$ belongs to $ C([0,1] : \mathcal{H}^1)$. As a consequence $u_0$ and $u_1$ are continuous on $\mathbb R$, so that the condition \eqref{cond:u0u1} makes sense.
\end{remark}

\begin{proof}[Proof of Theorem \ref{theo:S}]
Let 
\begin{equation}
    I_k(f) = \int_{\mathbb R}|x|^{2k}|f(x)|^2 dx.
\end{equation}
We will prove that $$I_k(u_0)<\infty \quad \text{and} \quad I_k(u_1)<\infty \quad \text{ for all } k\ge 1.$$ 
Then, by Corollary \ref{Cor. S}, it follows that $u\in C([0,1] : \mathcal S(\mathbb R)).$ 

To this end, we follow closely the argument in Section 5.1 in \cite{kulikov_uniqueness}.
Let $X_0>0$ be a large number to be fixed later. For $|x|\geq X_0$, the interval between two zeroes of $u_1$ containing $x$ has length at most $(2\lambda_1|x|^b)^{-1}$, where $0<b<\frac1\beta-1$ and  $\lambda_1=\lambda_1(X_0)$ can be chosen as large as we want by fixing $X_0$ large. 

Indeed, assume without loss of generality that $x>X_0$ and $x \in \big(n^{\beta},(n+1)^{\beta}\big)$. Then, by the mean value theorem we have
$$|(n+1)^\beta-n^\beta|\leq \beta n^{-(1-\beta)}.$$ Since $x\leq (n+1)^\beta\leq (2n)^\beta,$ we have that
$$n^{-(1-\beta)}\leq 2^{1-\beta}x^{-\frac{1-\beta}\beta} \leq 2^{1-\beta}x^{-b }x^{-\delta}\leq 2^{1-\beta}|x|^{-b}X_0^{-\delta},$$
where $\delta=\frac{1}{\beta}-1-b>0$.
Thus, $$|(n+1)^\beta-n^\beta|\leq \beta 2^{1-\beta}X_0^{-\delta}|x|^{-b} = (2\lambda_1|x|^{b})^{-1},$$ for $\lambda_1=\frac{2^{\beta}X_0^{\delta}}{4\beta}$. 

We argue exactly as in Section 5.1 in \cite{kulikov_uniqueness}  for $f=\tilde{u}_1$. Although our definitions of $a$ and $b$ differ, this does not affect the argument. Then, we get from estimate (5.6) in \cite{kulikov_uniqueness} that
\begin{align}\label{KNS:bound1}
     \int_{\mathbb R} |\tilde{u}_1(x)|^2 \Phi_R(|x|^{2b}) dx & \le \frac 1{\sqrt{\lambda_1}} \int_{\mathbb{R} }|\widehat{\tilde{u}}_1(\xi)|^2 \Phi_R(\xi^2)d\xi
+c\|\tilde{u}_1\|_{L^2}^2,\end{align}
where $$\Phi_R(t) = \begin{cases}
    t^k & t\in [0,R] \\ R^{k}+k R^{k-1}(t-R) & t>R,
\end{cases}$$
for $k\ge1.$
Similarly, replacing the role of $u_1$ and $u_0$, and working with the function $f=\tilde{v}_{-,1}=e^{ix^2/4}v_1=e^{ix^2/4}u_0$, instead of $\tilde{u}_1$, and for $0<a<\frac1\alpha-1,$ $\lambda_0 = \frac{2^{\alpha}X_0^{\delta}}{4\alpha}$, for some $\delta>0,$ we obtain that 
\begin{align}\label{KNS:bound2}
     \int_{\mathbb R} |\tilde{v}_{-,1}(x)|^2 \Phi_R(|x|^{2a}) dx &\le \frac 1{\sqrt{\lambda_0}} \int_{\mathbb{R} }|\widehat{\tilde{v}}_{-,1}(\xi)|^2 \phi_R(\xi^2)
+c\|u_0\|_{L^2}^2.\end{align}
We will now prove by induction that the estimates \eqref{KNS:bound1} and \eqref{KNS:bound2} imply that $u_0,\, u_1 \in L^2(x^{2k}dx )$
for all $k\ge 1.$
Indeed, since $u_0\in \mathcal{H}^1$, the property holds for $k=1$ thanks to Proposition \ref{WP} with $k=1$. 
Now, assume that $x^ku_0 \in L^2$ for some $k\in \mathbb N$, $k\ge 1$. We observe from Proposition \ref{Prop2}, in particular \eqref{bound on Gamma}, that 
\begin{align*}
    \int_{\mathbb{R} }|\widehat{\tilde{u}}_1(\xi)|^2 \Phi_R(\xi^2)d\xi \le \int_{\mathbb{R}}|\widehat{\tilde{u}}_1(\xi)|^2 \xi^{2k} d\xi  \le \int_{\mathbb{R}} |\partial_x^k \tilde{u}_1|^2 dx = (2)^{-k} \int_{\mathbb{R}}|\Gamma^k_1 u_1|^2 dx \le c(k,V) \left(\|x^k u_0\|_{L^2}^2 + \|u_0\|_{L^2}^2 \right). 
\end{align*}
Arguing as in the proof of Proposition \ref{Prop2}, the last estimate is rigorously justified assuming only that $x^ku_0 \in L^2$. This estimate combined with \eqref{KNS:bound1} yields that
the bound
\begin{align}\label{A2}
     \int_{\mathbb R} |\tilde{u}_1(x)|^2 \Phi_R(|x|^{2b}) dx  &\le c(k,V) \left( \|x^k u_0\|_{L^2}^2 + \|u_0\|_{L^2}^2\right)\end{align}
     holds.
     By letting $R\to+\infty$ and using the monotone convergence theorem, it follows that
\begin{align}\label{u_1 bk}
     \int_{\mathbb R} |x|^{2bk} |\tilde{u}_1|^2 dx  &\le c(k,V) \left( \|x^k u_0\|_{L^2}^2 + \|u_0\|_{L^2}^2\right) ,\end{align}
which implies that $x^{bk}u_1 \in L^2$. 
Arguing similarly for $\tilde{v}_{-,1}$, and using \eqref{KNS:bound2} instead \eqref{KNS:bound1}, we deduce that
\begin{align}\label{B}
     \int_{\mathbb R} |x|^{2aj} |\tilde{v}_{-,1}|^2 dx  &\le c(j,V) \left( \|x^j u_1\|_{L^2}^2 + \|u_0\|_{L^2}^2\right)
\end{align}
for $j=j(k)=\floor{bk}$. 
This estimate combined with \eqref{u_1 bk}, (recall also that $\tilde{v}_{-,1}=e^{ix^2/4}u_0$) implies
that
\begin{align}\label{B2}
     \int_{\mathbb R} |x|^{2aj(k)} |u_0|^2 dx  &\le c(k,V) \left( \|x^{j(k)} u_1\|_{L^2}^2+  \|u_0\|_{L^2}^2\right) \le c(k,V)\left(\|x^k u_0\|_{L^2}^2 + \|u_0\|_{L^2}^2\right). \end{align}
Next, we assume $x^ku_1 \in L^2$. We deduce from a similar reasoning that   
\begin{equation}\label{A3}
\int_{\mathbb R} |x|^{2b\hat{j}(k)} |\tilde{u}_1|^2 dx  \le c(k,V) \left( \|x^{\hat{j}(k)} u_0\|_{L^2}^2+  \|u_0\|_{L^2}^2\right) \le c(k,V)\left(\|x^k u_1\|_{L^2}^2 + \|u_0\|_{L^2}^2\right) ,
\end{equation}
for $\hat{j}(k) =\floor{ak}$. 

By hypothesis, we know that $0<\alpha,\beta<1$ satisfy $\alpha+\beta<1$, which is equivalent to the condition 
$${\left(\frac{1-\alpha}{\alpha}  \right)\left(\frac{1-\beta}{\beta} \right) }>1.$$
Thus we can choose $a<\frac{1}{\alpha}-1$ and $b<\frac{1}{\beta}-1$ such that $ab>1$.

Since 
$aj(k) = a\floor{bk} \ge a(bk-1),$  there exists an integer $k_0$ such that
$$aj(k) \ge k+1, \quad \forall \, k\ge k_0 .$$ 
Then, it follows from \eqref{B2} that
$$\int_{\mathbb R}|x|^{2(k+1)}|u_0|^2 \le \int_{\mathbb R}|x|^{2aj(k)}|u_0|^2 + \|u_0\|_{L^2}^2 \le c(k,V)\left(\|x^k u_0\|_{L^2} + \|u_0\|_{L^2}\right).$$
Arguing similarly, there exists $k_1 \in \mathbb N$ such that  
\[ b\hat{j}(k)\ge k+1, \quad \forall \, k\ge k_1, \] so by \eqref{A3}
$$\int_{\mathbb R}|x|^{2(k+1)}|u_1|^2 \le c(k,V)\left(\|x^k u_1\|_{L^2} + \|u_0\|_{L^2}\right) . $$
For $k \le K_0:=\max\{k_0,k_1\}$, we have directly 
$$\int_{\mathbb R} |x|^{2(k+1)}|u_0|^2dx+\int_{\mathbb R} |x|^{2(k+1)}|u_1|^2dx \le \|u_0\|_{L^2}^2 + \int_{\mathbb{R}} |x|^{2(K_0+1)}|u_0|^2 dx+\int_{\mathbb{R}} |x|^{2(K_0+1)}|u_1|^2 dx  <\infty.$$
This concludes the induction process. 
Therefore, we deduce that $u_1, u_0\in L^2(x^{2k}dx)$ for all $k\ge 1$.  Finally we apply Corollary \ref{Cor. S} to conclude that $u\in C([0,1] : \mathcal{S}(\mathbb R)).$
\end{proof}

Now, we prove Theorem \ref{theo:S:NLS}.
 The proof for the nonlinear case follows the same argument as for the potential case, and by adapting estimate \eqref{ener:est} on $\|\Gamma^k_t u\|_{L^2}$ in the case of the nonlinear Schr\"odinger equation \eqref{NLS}, (see \cite{kehle_uniqueness_2022,hayashi_solutions_1986, hayashi_solutions_1987,hayashi_nonlinear_1988}).

\begin{proof}[Proof of Theorem \ref{theo:S:NLS}]
We follow the same strategy as in the potential case. The argument is identical until \eqref{KNS:bound2}. By Proposition \ref{LWP:NLS:bound}, we know that if $u_0\in \mathcal{H}^1$, then $u_1\in \mathcal{H}^1,$ so the induction hypothesis holds for $k=1.$ Instead of the estimate \eqref{bound on Gamma}, we use the bound
    \begin{equation}
        \sup_{t\in [0,1]}\|\Gamma^k_t u\|_{L^2}\le c(k, \|u_0\|_{\mathcal{H}^1}) \|x^k u_0\|.
    \end{equation}
    This estimate corresponds to $(4.12)$ in the proof of Lemma 4.1 in \cite{kehle_uniqueness_2022} and is also rigorously justified by only assuming $x^ku_0\in L^2$.

    Then, the rest of the argument follows as for the potential case, and the result follows by Corollary \ref{Cor. S:NLS}.
\end{proof}

\section{Proof of the unique continuation result in the case of the potential} \label{sec:UC_potential}
This section is dedicated to the proof Theorem \ref{theo:UC:potential}. Below, we give a short overview of the proof.

We have seen in Section \ref{sec:S} that, for   potential $V \in L^{\infty}( [0,1] : C^{\infty}_b(\mathbb R))$, the solutions $u$ of \eqref{SP} emanating from initial data in $\mathcal{H}^1$ are in Schwartz for all time $t \in [0,1]$ if one assumes that $u(\pm n^{\alpha},0)=u(\pm n^{\beta},1)=0$, for all $n\in \mathbb N$ and $\alpha,\beta\in (0,1)$ with $0<\alpha+\beta <1$. In this section, we impose stronger decay assumptions on $V$, which in each case will provide a subset of $(\alpha,\beta)$ for which the unique continuation result from the zeroes holds. By using Theorem \ref{theo:S}, we know that $u \in C([0,1] : \mathcal{S}(\mathbb R))$,  so that all the computations below are rigorous. 

We set 
$$I_k(f) = \int_{\mathbb R} |x|^{2k}|f(x)|^2dx.$$ 
Moreover, let
$$0<a<\frac{1-\alpha}{\alpha},\quad  0<b<\frac{1-\beta}{\beta}.$$
We will prove that 
\begin{equation}\label{eq:goal} I_{\frac{pk}{2}}(u_0)\le c^k (k!)^{p\delta_0}, \quad I_{\frac{qk}{2}}(u_1)\le c^k (k!)^{q\delta_1},\end{equation} for $\frac1p+\frac1q=1,$ where $\delta_0$ and $\delta_1$ depend on $a$ and $b$. 

With estimate \eqref{eq:goal} in hand, we deduce Gaussian decay for $u_0$ and $u_1$ as soon as $p\delta_0<1$ and $q\delta_1<1$. Then, we will use Morgan's uncertainty principle, Theorem \ref{theo:EKPV2011}, to conclude that $u \equiv 0$. The conditions $p\delta_0<1$ and $q\delta_1<1$ will give us the final restrictions on $\alpha$
and $\beta$ by optimizing the choice of $p$ and $q$.

\subsection{Decay estimates from the zeros.}
Here we follow closely the argument in Section 5.2 in \cite{kulikov_uniqueness}. As in the proof of Theorem \ref{theo:S}, we have that for $|x|\ge X_0,$ the interval interval between two zeroes of $u_1$ containing $x$ has length at most $(2\lambda_1|x|^b)^{-1}$ for $\lambda_1=\frac{2^{\beta}X_0^{\delta}}{4\beta}$, $0<b<\frac{1}{\beta}-1$ and $\delta=\frac{1}{\beta}-1-b$. 

We argue exactly as in Section 5.2 in \cite{kulikov_uniqueness}  for $f=\tilde{u}_1$. Although our definitions of $a$ and $b$ differ, this does not affect the argument. Then, we get from estimate (5.11) in \cite{kulikov_uniqueness}, by taking $X_0$ large enough, that for $k\ge k_0$, 
\begin{equation*}
    I_{bk}(u_1)=I_{bk}(\tilde{u}_1) =\int |x|^{2bk}|\tilde{u}_1(x)|^2 dx \le \frac{1}{2}\int_{\mathbb{R}} |\xi|^{2k}|\widehat{\tilde{u}}_1(\xi)|^2 d\xi + c_{\beta}^k k^{\frac{2b}{b+1}k}I_0(\tilde{u}_1).
\end{equation*}
By using Plancherel's identity, it follows that 
\begin{align}\label{I_k u_1} \nonumber
    \int_{\mathbb R}|x|^{2bk} |\tilde{u}_1(x)|^2 dx &\leq \frac12 (2\pi)^{-2k}\int_\mathbb R  {|\partial_x^{k}\tilde{u}}_1 (x)|^2 dx
+c_{\beta}^k k^{\frac{2b}{b+1}k} I_0(u_1) \\ & \le \frac12 2^{-2k} \|\Gamma^k_1 u_1\|_{L^2}^2 + c_\beta^k k^{\frac{2b}{b+1}k}I_0(u_1).
\end{align}
Similarly, replacing the role of $u_1$ and $u_0$ and $\alpha$ and $\beta$ we deduce that,
\begin{equation} \label{I_k u_0}
I_{ak}(u_0)=\int_{\mathbb R}|x|^{2ak} |\tilde{v}_{-,1}(x)|^2 dx \leq\frac12 2^{-2k} \|\Gamma^k_{-1}u_0\|_{L^2}^2
+c_{\alpha}^k k^{\frac{2a}{a+1}k} I_0(u_0), 
\end{equation}
for $0<a<\frac1\alpha -1$, and $k\ge k_0,$ by increasing $k_0$ if necessary.

\smallskip
\subsection{Controlling the evolution of $\Gamma_t^k u$ in $L^2$}
The aim of this subsection is to control the $L^2$-norms of $\Gamma_t^ku$ and $\Gamma_{-t}^k\tilde{v}_-$ that appear on the right-hand side of \eqref{I_k u_1}-\eqref{I_k u_0}.

\begin{proposition} \label{lemma:EE:Gamma}
For any integer $k\ge 1$
\begin{equation}\label{Gamma estimate}
       \frac12\sup_{t\in [0,1]}\|\Gamma^k_tu(t)\|_{L^2}\leq \|x^ku_0\|_{L^2} + c^{k}(k!)^{\tau}\|u_0\|_{L^2}, 
\end{equation}
\begin{equation}\label{Gamma estimate for v}
       \frac12\sup_{t\in [0,1]}\|\Gamma^k_{-t}v(t)\|_{L^2}\leq \|x^ku_1\|_{L^2} + c^{k}(k!)^\tau\|u_0\|_{L^2}. 
\end{equation}
where $$\tau = \begin{cases}
1 , & \text{when} \ 0 \leq \sigma < 1,\\
    \sigma^+, & \text{when}  \ \sigma \geq 1, 
\end{cases}$$
where $\sigma^+$ is any number slightly larger than $\sigma$.
\end{proposition}
\medskip


\begin{proof}[Proof of Proposition \ref{lemma:EE:Gamma}]
To prove \eqref{Gamma estimate}, we follow the same strategy as in Section \ref{sec:S}, but now we will also track the constant to get a more precise decay estimate. By using \eqref{ener:est} and H\"older's inequality, we get
\begin{align*}
    \frac{1}{2}\frac{d}{dt}\|\Gamma^k_t u(t)\|_{L^2}^2 &\leq \sum_{l=0}^{k-1}\binom{k}{l}(2t)^{2k} \|\partial_x^{k-l}V\|_{L^\infty} \| \partial_x^l \tilde{u}(t)\|_{L^2}\|\partial_x^k \tilde{u}(t)\|_{L^2} \\ & = \sum_{l=0}^{k-1}\binom{k}{l}(2t)^{k} \|\partial_x^{k-l}V\|_{L^\infty}\|\partial_x^{l} \tilde{u}(t)\|_{L^2}\|\Gamma^k_t u(t)\|_{L^2}.
    \end{align*}
This implies that
\begin{equation*}
    \frac{d}{dt}\|\Gamma^k_tu(t)\|_{L^2}\leq \sum_{l=0}^{k-1}\binom{k}{l}(2t)^{k}\|\partial_x^{k-l}V\|_{L^\infty}\|\partial_x^l \tilde{u}(t)\|_{L^2}.
\end{equation*}
For all $1\le l \le k-1, $ let $\theta=\theta_l>0$ to be fixed later. Then 
\begin{align*}
    \|\partial_x^l \tilde{u}\|_{L^2} = (2\pi)^{l}\|\xi^l\widehat {\tilde{u}}(t)\|_{L^2}&=(2\pi)^l \|{\bf 1}_{|\xi| \le \theta} \xi^l \widehat{\tilde{u}}(t)\|_{L^2} + (2\pi)^l\|{\bf 1}_{|\xi| > \theta} \xi^l \widehat{\tilde{u}}(t)\|_{L^2} \\ &\le (2\pi)^l\theta^{l}\|u_0\|_{L^2} + (2\pi)^{l-k} \theta^{-(k-l)}\|\partial_x^k \tilde{u}(t)\|_{L^2}\\ & = (2\pi \theta)^l\|u_0\|_{L^2} + (2\pi \theta)^{-(k-l)}(2t)^{-k}\|\Gamma_t^k u(t)\|_{L^2}.
\end{align*}
By using that $\sup_{t\in [0,1]}\|\partial_x^{l}V\|_{L^\infty}\le A^{l+1} (l!)^{\sigma},$ it follows that
\begin{align}\label{ddt Gamma}
    \frac{d}{dt}\|\Gamma^k_tu(t)\|_{L^2}&\leq  A\sum_{l=0}^{k-1}\binom{k}{l}(2t)^{k}A^{k-l} \left((k-l)!\right)^{\sigma} (2\pi\theta)^l\|u_0\|_{L^2} \nonumber \\ & \quad  + A\sum_{l=1}^{k-1}\binom{k}{l}(2\pi A^{-1}\theta)^{-(k-l)}\left((k-l)!\right)^{\sigma}\|\Gamma_t^k u(t)\|_{L^2},
\end{align} and integrating \eqref{ddt Gamma} from $0$ to $t$, we get that
\begin{align} \label{rf}
   \sup_{t\in [0,1]} \|\Gamma^k_tu(t)\|_{L^2}&\leq \|x^ku_0\|_{L^2} + A\sum_{l=0}^{k-1}\binom{k}{l}2^{k}A^{k-l} \left((k-l)!\right)^{\sigma}(2\pi \theta)^l\|u_0\|_{L^2} \nonumber \\ & \quad + A\sum_{l=1}^{k-1}\binom{k}{l}(2\pi A^{-1}\theta)^{-(k-l)}\left((k-l)!\right)^{\sigma}\sup_{t\in [0,1]}\|\Gamma_t^k u(t)\|_{L^2}.
\end{align}

We let $J_0\in \mathbb N$, $J_0=J_0(\sigma) \ge 2$, and $\rho>0$ to be fixed later. Let $$\theta=\theta_l:=\begin{cases}
    k^\sigma \frac{A}{2\pi}\rho^{-1}, & 1\le l\le \floor{\frac k {J_0}} \\
    l^{\max\{1,\sigma\}} \frac{A}{2\pi}\rho^{-1}, & \floor{\frac k {J_0}}+1 \le l \le k-1.
\end{cases}$$
We will estimate the two sums separately. By Stirling's formula it follows, in the case where $1 \le l \le \floor{\frac{k}{J_0}}$,
\begin{equation} \label{est:II.1} A\sum_{l=0}^{\floor{\frac k{J_0}}}\binom{k}{l}2^{k}A^{k-l} \left((k-l)!\right)^{\sigma}(2\pi \theta)^l \le  Ak! (2 A)^k \sum_{l=0}^{\floor{\frac k{J_0}}} \left((k-l)!\right)^{\sigma-1} \rho^{-l} k^{\frac{k\sigma}{J_0}} \le \begin{cases}
    c^k(k!)^{\sigma+\frac 1 {J_0}}, & \sigma\le 1, \\
    c^k (k!)^{(1+\frac1{J_0})\sigma}, & \sigma >1.
\end{cases}
\end{equation}
In the case where $ l> \floor{\frac k{J_0}}$, we need to bound
  \begin{equation*} 
  A\sum_{l=\floor{k/J_0}+1}^{k-1}\binom{k}{l}2^{k}A^{k-l} \left((k-l)!\right)^{\sigma}(2\pi \theta)^l \le  c^k \sum_{l=\floor{\frac k{J_0}}+1}^{k-1} k!\left((k-l)!\right)^{\sigma-1} l^{(\max\{1,\sigma\}-1)l} .
\end{equation*}
When $\sigma \le 1$, we get the rough bound
\begin{equation} \label{est:II.2} A\sum_{l=\floor{k/J_0}+1}^{k-1}\binom{k}{l}2^{k}A^{k-l} \left((k-l)!\right)^{\sigma}(2\pi \theta)^l \le  c^k k! .
\end{equation}
When $\sigma>1$, we observe from the Stirling formula that
$$(k-l)!l^l \le (k-l)^{(k-l)}l^l= \frac{k^k}{k^k}(k-l)^{(k-l)}l^l = 
\left(\frac{(k-l)}{k}\right)^{k-l} \left(\frac{l}{k}\right)^{l} k^k\le c^k k! .$$
Thus, 
\begin{equation} \label{est:II.3} A \sum_{l=\floor{k/J_0}+1}^{k-1}\binom{k}{l}2^{k}A^{k-l} \left((k-l)!\right)^{\sigma}(2\pi\theta)^l \le  c^k k! \sum_{l=\floor{\frac k{J_0}}+1}^{k-1} \left((k-l)! l^l\right)^{\sigma-1} \le c^k (k!)^\sigma .
\end{equation}

For the other sum in \eqref{rf}, we will prove the following claim. 
\begin{claim} \label{claim:S} \textit{For $l, k \in \mathbb N$ with $1 \le l \le k-1$, let $$S_{l,k}=\binom{k}{l}\left((k-l)!\right)^{\sigma}(2\pi A^{-1}\theta_l)^{-(k-l)} =\begin{cases}
   \binom{k}{l}\left((k-l)!\right)^{\sigma}k^{-\sigma(k-l)}\rho^{k-l}, & 1\le l \le \floor{\frac k {J_0}}, \\
    \binom{k}{l}\left((k-l)!\right)^{\sigma}l^{-\max\{1,\sigma\}(k-l)}\rho^{k-l}, & \floor{\frac{k}{J_0}}+1\le l \le k-1.
\end{cases}$$
Then, there exists $\rho_0=\rho_0(\sigma,A)>0$ such that, for all $0<\rho<\rho_0$, for all $k \ge 2$,
\begin{equation} \label{bound:S}
\sum_{l=1}^{k-1}S_{l,k} \le \frac1{2A} .
\end{equation}}
\end{claim}
Assuming the claim, it follows that
\begin{align*}
A\sum_{l=1}^{k-1}\binom{k}{l}(2\pi A^{-1}\theta)^{-(k-l)}\left((k-l)!\right)^{\sigma}\sup_{t\in [0,1]}\|\Gamma_t^k u(t)\|_{L^2} \le \frac12\sup_{t\in [0,1]}\|\Gamma^k_t u(t)\|_{L^2},
\end{align*}
and combining this with \eqref{est:II.1} and \eqref{est:II.2}, we conclude the proof of \eqref{Gamma estimate}. 
It remains to prove Claim \ref{claim:S}.
\begin{proof}[Proof of Claim \ref{claim:S}]
First let $M\in \mathbb N, M>J_0.$ Let $m=k-l$. Observe that for $k\le M$, and $1 \le l \le k-1$, since $( A(2\pi\theta)^{-1})^m\le \rho^m$ as long as $\sigma \ge 0$, we always have
$$S_{l,k}\le\binom{k}{m}(m!)^\sigma \rho^m \le 2^M (M!)^{\sigma} \rho^m \le \left(2^M (M!)^{\sigma} \rho\right)^m < \left(\frac 1{4A}\right)^m $$ for $\rho<\frac 1{4A} \left(2^M (M!)^{\sigma} \right)^{-1}.$
Thus we can assume $k\ge M>J_0.$

We split the proof into two cases depending on whether $l$ is larger or smaller than $\floor{\frac{k}{J_0}}$.

\medskip

\noindent \emph{Case} $ \floor{\frac k {J_0}} \leq l\leq k-1.$ Then, $l\ge \frac{k}{J_0}-1\ge \frac k {J_0} -\frac k M = k\left( \frac{M-J_0}{MJ_0} \right) =:kN_0.$ (recall that $k\ge M$). Let $m=k-l.$ Since $\binom{k}{l}=\binom{k}{m}\leq \frac{k^m}{m!}$, we have
\begin{align*}
    S_{l,k} & = \binom{k}{m} \left(m!\right)^\sigma l^{-\max\{1,\sigma\}m}\rho^{m}\leq \frac{k^m}{m!}(m!)^{\sigma}\left(\frac{\rho}{l^{\max\{1,\sigma\}}}\right)^{m} .
\end{align*}
Then, when $\sigma\le 1$, it follows that 
\begin{align*}
    S_{l,k} \leq k^m \left(\frac{\rho}{kN_0}\right)^{m} = (N_0^{-1}\rho)^m<\left(\frac1{4A}\right)^m,
\end{align*}
for $\rho<\frac{N_0}{4A}.$
When $\sigma>1,$ 
$$S_{l,k} \le k^{m(1-\sigma)}(m!)^{\sigma-1}(N_0^{-\sigma} \rho)^m.$$
Since $m! \le m^m,$ it follows that
$$S_{l,k}= k^{m(1-\sigma)}m^{m(\sigma-1)}(N_0^{-\sigma} \rho)^m \le e^{m(1-\sigma)(\log{k} -\log{m}) }(N_0^{-\sigma} \rho)^m \le (N_0^{-\sigma} \rho)^m<\left( \frac1{4A}\right)^m$$ for $\rho <\frac{N_0^{\sigma}}{4A}.$

Thus, for $\floor{\frac k{J_0}} < l \le k-1,$
$S_{l,k}\le (\frac1{4A})^{m}$ for $\rho <\rho_1:=\frac1{4A} \min\{ N_0, N_0^\sigma, (2^M(M!)^\sigma)^{-1}\}.$

\medskip

\noindent \emph{Case} $1\leq l \leq \floor{\frac k{J_0}}.$ Then $k\geq J_0 l$. First, we observe that $S_{l,k}\geq S_{l,{k+1}}$.  Indeed,
\begin{align*}
  \frac{S_{l,k+1}}{S_{l,k}} = \frac{\binom{k+1}{l} \left((k+1-l)!\right)^{\sigma}(k+1)^{-\sigma(k+1-l)}\rho^{k+1-l}}{\binom{k}{l}\left((k-l)!\right)^{\sigma}k^{-\sigma(k-l)}\rho^{k-l}}=\frac{k+1}{k+1-l} (k+1-l)^{\sigma}\left(\frac{k}{k+1}\right)^{\sigma(k-l)}(k+1)^{-\sigma}\rho.
\end{align*}
When we fix $l\leq k/J_0$, the function $: k\mapsto \frac{k+1}{k+1-l}$ is decreasing, so $$\frac{k+1}{k+1-l}\leq \frac{J_0l+1}{J_0l+1-l}=\frac{J_0l+1}{l(J_0-1) +1} = 1+\frac l{(J_0-1)l + 1} \le 1+\frac 1{J_0-1}\le 2,$$ for $J_0\ge 2$.
Thus,
\begin{equation*}
   \frac{S_{l,k+1}}{S_{l,k}}\leq 2 \frac{(k+1-l)^{\sigma}}{(k+1)^\sigma}\left(\frac{k}{k+1}\right)^{\sigma(k-l)}\rho \le 2\rho,
\end{equation*}
since $\sigma \ge 0$.
So, 
$ \frac{S_{l,k+1}}{S_{l,k}}< 1 $ when $\rho <\frac12.$

It follows that $S_{l,k}\leq S_{l,J_0l}$ for $k \ge J_0l$.
Moreover, we have by using the Stirling formula,
\begin{align*}S_{l,J_0l}&=\binom{J_0l}{l}(((J_0-1)l)!)^\sigma (J_0l)^{-\sigma(J_0-1) l}\rho^{(J_0-1)l} \\ &\le \frac{(J_0 l)^l}{l!}  (((J_0-1)l)!)^\sigma (J_0l)^{-\sigma(J_0-1) l}\rho^{(J_0-1)l} \\ & \le c^l J_0^l (((J_0-1)l))^{\sigma(J_0-1)l} (J_0l)^{-\sigma(J_0-1) l}\rho^{(J_0-1)l} \\
& \le c^{l}J_0^l \rho^{(J_0-1)l} \le \left(\frac1{4A}\right)^l
\end{align*} 
for $\rho<\rho_2:=(4c AJ_0)^{-(J_0-1)}.$
This concludes the proof of the claim by choosing \[0<\rho<\rho_0=\min\{\rho_1,\rho_2\} .\]
\end{proof}
Therefore, we deduce gathering \eqref{est:II.1}, \eqref{est:II.2}, \eqref{est:II.3} and Claim \ref{claim:S} that
\begin{equation}
       \frac12\sup_{t\in [0,1]}\|\Gamma^k_tu(t)\|_{L^2}\leq \|x^ku_0\|_{L^2} + c^{k}(k!)^{\widetilde{\tau}}\|u_0\|_{L^2}, 
\end{equation}
where 
$$\widetilde{\tau} = \begin{cases}
    \max\{\sigma+\frac 1 {J_0},1\}, & 0 \le \sigma < 1, \\ \sigma (1+\frac 1{J_0}), & \sigma  \geq 1.
\end{cases}$$
When $0 \leq \sigma < 1,$ we choose $J_0$ large enough so that $\sigma+\frac 1 {J_0}<1$.
When $\sigma \geq 1$,  we choose $J_0$ large enough so that $\sigma\big(1+\frac 1 {J_0}\big)<\sigma^+$, for any number $\sigma^+>\sigma$. This concludes the proof of \eqref{Gamma estimate}. Arguing similarly, for $v(x,t) = u(x,1-t)$ and $\Gamma^k_{-t}v(t),$ we deduce \eqref{Gamma estimate for v}. This concludes the proof of Proposition \ref{lemma:EE:Gamma}.

\bigskip

\subsection{Conclusion of the proof}
By inserting \eqref{Gamma estimate} and \eqref{Gamma estimate for v} into \eqref{I_k u_0} and \eqref{I_k u_1} we obtain, for $k\ge k_0$
\begin{align}\label{Ik:u1:b}
    \int_{\mathbb{R}} |x|^{2bk}|u_1(x)|^2 dx & \le \frac{1}{2}\|x^k u_0\|_{L^2}^2 + \left(c^k(k!)^{2\tau} + c_\beta^k k^{\frac{2b}{b+1} k}\right) I_0(u_0)
\end{align}
and 
\begin{align}\label{Ik:u0:a}
    \int_{\mathbb{R}} |x|^{2ak}|u_0(x)|^2 dx & \le \frac{1}{2}\|x^k u_1\|_{L^2}^2 + \left(c^k(k!)^{2\tau} + c_\alpha^k k^{\frac{2a}{a+1} k}\right) I_0(u_0).
\end{align}

Let $1<p,q<\infty$ with $\frac1p+\frac1q=1$ to be fixed later, and $m\in \mathbb N$. Applying \eqref{Ik:u1:b} with $k=\floor{\frac{qm}{2b}}$ and using Stirling's formula imply
\begin{equation}\label{u_1 pq}
    \int_{\mathbb{R}} |x|^{qm}|u_1(x)|^2 dx \leq  \frac 12 \int |x|^{\frac{qm}{b}}|u_0(x)|^2 dx  + \left(c^m(m!)^{\frac{q\tau}{b}} + c_\beta^m m^{\frac{q}{b+1} m}\right) I_0(u_0).
\end{equation} 
Applying \eqref{Ik:u0:a} with $k=\floor{\frac{qm}{2ab}}$ and using Stirling's formula imply
\begin{equation}\label{u_0 pq}
      \int_{\mathbb{R}} |x|^{\frac{qm}{b}}|u_0(x)|^2 dx \le \frac{1}{2}\int_{\mathbb R}|x|^{\frac{qm}{ab}}|u_1(x)|^2 + \left(c^m(m!)^{\frac{q\tau}{ab}} + c_\alpha^m m^{\frac{q}{(a+1)b} m}\right) I_0(u_0).
\end{equation}
Inserting \eqref{u_0 pq} into \eqref{u_1 pq} yields
\begin{equation}
    \int_{\mathbb R} |x|^{qm}|u_1(x)|^2dx \le \frac{1}{4}\int|x|^{\frac{qm}{ab}} |u_1(x)|^2 + c^m\left( (m!)^{\frac{q\tau}{ab}} + m^{\frac{q}{(a+1)b}m}+(m)!^{\frac{q\tau}{b}} + m^{\frac{q}{b+1}m} \right).
\end{equation}
Recall that condition $\alpha+\beta<1$ is equivalent to
$$\left(\frac{1-\alpha}{\alpha}\right)\left(\frac{1-\beta}{\beta}\right)>1.$$ Therefore, it is possible to choose $a<\frac{1-\alpha}{\alpha}$ and $b<\frac{1-\beta}{\beta}$ such that $ab>1.$ Then 
$$\frac{1}{4}\int|x|^{\frac{ql}{ab}} |u_1(x)|^2 \le \frac{1}{4}\int |x|^{ql}|u_1(c)|^2 dx + I_0(u_0).$$ 
It follows from Stirling's formula that, for $k \ge k_0$,
\begin{equation}\label{final:decay:1}
    \int_{\mathbb R} |x|^{qk}|u_1(x)|^2dx \le e^{q\delta_1 k\log{k}+O(k)},
\end{equation}
where 
$\delta_1=\max\left\{\frac{\tau}{ab}, \frac{1}{(a+1)b}, \frac{\tau}{b},\frac{1}{b+1} \right\}=\max\left\{\frac{\tau}{ab}, \frac{\tau}{b} \right\}$. Here, we used that $\tau \ge 1$ and that $a,b>0.$

Applying \eqref{Ik:u0:a} with $k=\floor{\frac{pm}{2a}}$ and \eqref{Ik:u1:b} with $k=\floor{\frac{pm}{2ab}}$ and arguing  similarly, for all $k \ge k_0$,
\begin{equation}\label{final:decay:0}
    \int_{\mathbb R} |x|^{pk}|u_0(x)|^2dx \le e^{p\delta_0 k\log{k}+O(k)},
\end{equation}
where 
$\delta_0=\max\left\{\frac{\tau}{ab}, \frac{1}{(b+1)a}, \frac{\tau}{a},\frac{1}{a+1} \right\}=\max\left\{\frac{\tau}{ab}, \frac{\tau}{a}\right\}$.

Note that estimates \eqref{final:decay:0} and \eqref{final:decay:1} actually hold for all $k\in \mathbb N$. Indeed, in the case where $k \le k_0$, we use the fact that $u \in C([0,1] : \mathcal{S}(\mathbb R))$ (see Corollary \ref{Cor. S}), so that $I_k(u_0) +I_k(u_1)\le c(k_0,\|u_0\|_{\mathcal{H}^{k_0}})$.

Then, we conclude that
\begin{equation} \label{est:gauss:u1}
\int_{\mathbb R} e^{c|x|^q}|u_1(x)|^2dx =\sum_{k=0}^{\infty} \frac{c^k}{k!} \int|x|^{qk}|u_1(x)|^2 dx \le \sum_{k=0}^{\infty} e^{(q\delta_1-1)k\log{k} + O(k)} <+\infty 
\end{equation}
and 
\begin{equation} \label{est:gauss:u0}
\int_{\mathbb R} e^{c|x|^p}|u_0(x)|^2dx =\sum_{k=0}^{\infty} \frac{c^k}{k!} \int|x|^{pk}|u_0(x)|^2 dx \le \sum_{k=0}^{\infty}e^{(p\delta_0-1)k\log{k} + O(k)} <+\infty 
\end{equation}
if we require that $q\delta_1<1$ and $p\delta_0<1.$

Now, we optimize the choice of $p$ and $q$ to get the best conditions for $\alpha,\beta$. We choose $p$ and $q$ such that
$$p\delta_0=q\delta_1 \iff \frac{q}{q-1} \delta_0=q\delta_1 \iff q=\frac{\delta_0+\delta_1} {\delta_1}.$$ 
Then, it follows that  
\begin{equation} \label{cond:delta}
q\delta_1<1 \ \text{and} \ p\delta_0 <1 \iff \delta_0+\delta_1 <1.
\end{equation}

Below, we look at the meaning of condition \eqref{cond:delta} for $\alpha$ and $\beta$.

If $\delta_0=\frac{\tau}{ab}$, then we must have $$\frac{\tau}{ab}\ge\frac{\tau}{a} \iff b\le1.$$
But $b\le 1$ implies $\delta_1 \ge \frac{\tau}{b}\ge 1$,
so that $\delta_1+\delta_0 \ge 1$, which contradicts \eqref{cond:delta}. 

Similarly, the case $\delta_1 =\frac{\tau}{ab} \iff a \le 1$ is impossible. 

Thus, in order for \eqref{cond:delta} to be satisfied, we must have $\delta_1=\frac{\tau}{b}$ and $\delta_0=\frac{\tau}{a}$, with $a>1$ and $b>1$. 
In this case, 
\begin{equation} \label{cond:delta.2}
    \delta_0+\delta_1 <1 \iff \tau \left(\frac{1    }{a}+\frac{1}{b}\right)<1.
\end{equation}
First, observe that 
\begin{equation*}
\frac{\alpha}{1-\alpha}<\frac1a<1 \quad \text{and} \quad \frac{\beta}{1-\beta}<\frac1b<1,
\end{equation*}
so we must have $\alpha,\beta \in (0,\frac{1}{2})$.
Now, assume that
\begin{equation} \label{cond:alpha-beta-tau}
 \tau \left(\frac{\alpha    }{1-\alpha}+\frac{\beta}{1-\beta}\right)<1.
\end{equation}
Since $\alpha< \frac12$ and $\beta<\frac12$, then both $\frac{1-\alpha}\alpha>1$ and $\frac{1-\beta}{\beta}>1$. Hence, we may choose $1<a<\frac{1-\alpha}{\alpha}$ and $1<b<\frac{1-\beta}{\beta}$ such that \eqref{cond:delta.2} is satisfied, which in turn implies \eqref{cond:delta}.

Finally, we conclude  by using \eqref{est:gauss:u1}-\eqref{est:gauss:u0} and applying the Morgan principle proved by Escauriaza, Kenig, Ponce and Vega (see Theorem \ref{theo:EKPV2011}) that $u \equiv 0$. This finishes the proof of Theorem \ref{theo:UC:potential}.

\end{proof}

\section{Proof of the unique continuation result in the nonlinear case} \label{sec:UC_NLS}
We dedicate this section to the proof of Theorem \ref{thm:main}.  Without loss of generality, we assume that $\lambda=1$ in \eqref{NLS}. Below, we give a brief overview of the proof. 

From Theorem \ref{theo:S:NLS} in Section \ref{sec:S}, we know that solutions of \eqref{NLS} with $u_0\in \mathcal{H}^1$ that satisfies \eqref{zero_cond} for $0<\alpha + \beta <1$ are in Schwartz for all time $t\in [0,1].$ The proof follows the same overall strategy as in the preceding section. However, instead of controlling the quantities $I_k(u_0)$ and $I_k(u_1)$, we control the quantities
\begin{equation}
    E_k(f) :=I_k(f) + I_k(\partial_x f)=\int |x|{^{2k}}\left(|f(x)|^2 + |\partial_x f(x)|^2\right)dx.
\end{equation}
Note that, for all  $0\le k\le  r$,
\begin{equation}\label{Almost:monotonicity}
E_k(f) \le E_r(f) + E_0(f).
\end{equation}
It is worth noting that $E_0(f) = \|f\|_{H^1}^2 \le \|f\|_{\mathcal{H}^1}^2$, and that
\begin{equation} \label{est:xjf}
    \|x^kf\|_{H^1}^2 = \|x^k f\|_{L^2}^2 + \|\partial_x (x^k f) \|_{L^2}^2 \le 2(k^2+1) E_k (f) + 2k^2 \|f\|_{L^2}^2.
\end{equation}

We will again use the functions $\tilde{u}_1=e^{-ix^2/4}u_1$ and $\tilde{v}_{-,1}=e^{ix^2/4}u_0$. Differentiating the exponential, we obtain that 
\begin{equation} \label{comparing:with:tilde1}
    E_k(u_1) \le c \left(E_{k+1}(\tilde{u}_1) + E_0(\tilde{u}_1) \right) \le c \left(E_{k+1}(\tilde{u}_1) + \|u_1\|_{\mathcal{H}^1}^{2} \right),
\end{equation}
and
\begin{equation}\label{comparing:with:tilde:0}
    E_k(u_0) \le c \left(E_{k+1}(\tilde{v}_{-,1}) + E_0(\tilde{v}_{-,1}) \right) \le c \left(E_{k+1}(\tilde{v}_{-,1}) + \|u_0\|_{\mathcal{H}^1}^2 \right),
\end{equation}

Let $\alpha, \beta \in (0,\frac12)$. We work with $a$, $b$ satisfying 
\[ 1<a<\frac{1-\alpha}{\alpha}, \quad 1<b<\frac{1-\beta}{\beta} .\]

The main goal is to prove that 
\begin{equation} \label{E_pqk:u_0:u_1} E_{\frac{pk}2}(u_0) \le c^{k} (k!)^{p\delta_0}, \quad E_{\frac{qk}{2}}(u_1) \le c^k (k!)^{q\delta_1},\end{equation}
for $\frac1p+\frac1q=1,$ and $\delta_0=\frac1a$, $\delta_1=\frac1b$.  
These estimates provide Gaussian decay for $u_0$ and $u_1$ as soon as $p\delta_0<1$ and $q \delta_1<1$. We will then optimize $p$ and $q$ as in the proof of Theorem \ref{theo:UC:potential} and conclude by using the Morgan uncertainty principle (Theorem \ref{theo:EKPV2011}).

\subsection{Decay estimates from the zeros}\label{sec:zero-gap}

\begin{proposition}\label{prop:zero-gap}
Let $0<\alpha,\beta<1$ and $1<b<\frac{1}{\beta}-1$, $1<a<\frac{1}{\alpha}-1$. Then, there exist $c\ge 1$ independent of $u$ and $k_0 \in \mathbb N$ with $k_0 \ge 2$, such that for $k \ge k_0$,
\begin{equation}\label{eq:zero-gap-H1:u_1}
 E_{bk}(\tilde{u}_1)
 \le \frac12\left\|\tilde{u}_1\right\|_{H^{k+1}}^2
 +c^kk^{\frac{2b}{b+1}k}\|u_1\|_{\mathcal{H}^1}^2.
\end{equation}
\begin{equation}\label{eq:zero-gap-H1:u_0}
 E_{ak}(\tilde{v}_{-,1})
 \le \frac12\left\|\tilde{v}_{-,1}\right\|_{H^{k+1}}^2
 +c^kk^{\frac{2a}{a+1}k} \|u_0\|_{\mathcal{H}^1}^2.
 \end{equation}
\end{proposition}

\begin{proof}
We begin with the estimate from \cite[(5.11)]{kulikov_uniqueness}, used earlier. For $|x|\ge X_0$ and $x \in (n^\beta, (n+2)^\beta)$, we have by the mean value theorem that
\begin{equation*}
 (n+2)^\beta-n^\beta
 \le 2\beta n^{\beta-1}
 \le 2\beta 3^{1-\beta}|x|^{1-1/\beta}
 \le (2\lambda_1|x|^b)^{-1},
\end{equation*}
for $\lambda_1>1$. Hence, the consecutive zeros of $\tilde{u}_1$ satisfy
the gap hypothesis of \cite[(5.11)]{kulikov_uniqueness}. In particular, there exists $k_0 \in  \mathbb N$ with $k_0 \ge 2$ such that, for all $k \ge k_0$,
\begin{equation} \label{Ebk:1}
    I_{bk}(\tilde{u}_1) \le \frac{1}{2}\int |\partial_x^k \tilde{u}_1 (x)|^2 dx + c_\beta^{k}k^{\frac{2b}{b+1}k}\|u_1\|_{L^2}^2 .
\end{equation}
We now apply the same estimate to the derivative $\partial_x \tilde{u}_1$. As in \cite[Section~2.1]{ramos_fourier_2022}, Rolle's theorem gives, for
each $n\ge1$,
\begin{align*}
 \{x\in(n^\beta,(n+1)^\beta):
       \partial_x\operatorname{Re}\tilde{u}_1(x)=0\}&\ne\varnothing,\\
 \{x\in(n^\beta,(n+1)^\beta):
       \partial_x\operatorname{Im}\tilde{u}_1(x)=0\}&\ne\varnothing.
\end{align*}
The same argument applies between consecutive negative zeros. Consecutive zeros chosen
in these intervals lie within two consecutive gaps of $\{\pm n^\beta:n\ge1\}$,
so their distance is at most $(n+2)^\beta-n^\beta$ on the positive
half-line, with the corresponding bound on the negative half-line. The
previous gap estimate then yields
\begin{align}
\label{Ebk:2} I_{bk}(\partial_x\operatorname{Re}f) &\le\frac12\int |\partial_x^{k+1}\operatorname{Re}\tilde{u}_1(x)|^2
 +c^{k}k^{\frac{2b}{b+1}k} \left\|\partial_x\operatorname{Re}\tilde{u}_1\right\|_{L^2}^2,\\
 \label{Ebk:3}I_{bk}(\partial_x\operatorname{Im}\tilde{u}_1)
 &\le\frac12\int |\partial_x^{k+1}\operatorname{Im}\tilde{u}_1(x)|^2 dx
 +c^{k}k^{\frac{2b}{b+1}k}
\left\|\partial_x\operatorname{Im}f\right\|_2^2.
\end{align}
Combining the estimates \eqref{Ebk:1}, \eqref{Ebk:2} and \eqref{Ebk:3}, we deduce that
\begin{align*}
    E_{bk}(\tilde{u}_1)& \le \frac{1}{2} \left( \|\partial_x^k \tilde{u}_1\|_{L^2}^2 +\|\partial_{x}^{k+1} \tilde{u}_1\|_{L^2}^2\right) + c^k k^{\frac{2b}{b+1}k}\|\tilde{u}_1\|_{H^1}^2 \\ &\le \frac12\|\tilde{u}_1\|_{H^{k+1}}^2 + c^k k^{\frac{2b}{b+1}k}\|u_1\|_{\mathcal{H}^1}^2,
\end{align*}
which concludes the proof of \eqref{eq:zero-gap-H1:u_1}. The proof for \eqref{eq:zero-gap-H1:u_0} follows by the same argument for $u_0$ instead of $u_1.$
\end{proof}

\subsection{Control of the evolution of $\Gamma^k_t$ in $L^2$}\label{sec:transport}

In order to control the quantities $\left\|\tilde{u}_1\right\|_{H^{k+1}}$ and $\left\|\tilde{v}_{-,1}\right\|_{H^{k+1}}$ in \eqref{eq:zero-gap-H1:u_1}-\eqref{eq:zero-gap-H1:u_0},  we estimate the $L^2$-norm of the vector field operators $\Gamma^k_t$ in $L^2$ along the flow. Recall that $\Gamma_t$ is defined in \eqref{Gamma:form.1}; in particular, $\Gamma_0=x$. 

\begin{proposition} \label{prop:indu:Gammaj}
There exists $K=K(\|u_0\|_{\mathcal{H}^2}) \ge 1$ such that the following is true.

 (i) Assume that there exist $B \ge 1$ and $J \in \mathbb N$ with $J \ge 1$ such that
\begin{equation}\label{ass:Ej}
 E_j(u_0)\le (B^j j!)^2,\qquad 1\le j\le J .
\end{equation}
Then,
\begin{equation}\label{est:Gammaj_nl}
 \sup_{0\le t\le1}\left\|\Gamma_t^ju(t)\right\|_{H^1}
 \le (KB)^j j!,\qquad 1\le j\le J.
\end{equation}

(ii) Assume that there exist $B \ge 1$ and $J \in \mathbb N$ with $J \ge 1$ such that
\begin{equation}\label{ass:Ej1}
 E_j(u_1)\le (B^j j!)^2,\qquad 1\le j\le J .
\end{equation}
Then,
\begin{equation}\label{est:Gammaj_nl-}
 \sup_{0\le t\le1}\left\|\Gamma_{-t}^ju(1-t)\right\|_{H^1}
 \le (KB)^j j!,\qquad 1\le j\le J.
\end{equation}

\end{proposition}

In order to prove Proposition \ref{prop:indu:Gammaj}, we need the following technical lemma. 

\begin{lemma} \label{lemma:tec:convol}
There exists a universal constant $c>0$ such that for any $j\in \mathbb N$, $j \ge 2$, we have
\begin{align}
 \sum_{l=1}^{j-1}\frac{1}{(l+1)^2(j-l+1)^2}
 &\le\frac{c}{(j+1)^2},\label{lemma:tec:convol.1}\\
 \sum_{\substack{j_1+j_2+j_3=j\\j_1,j_2,j_3\ge1}}
 \frac{1}{(j_1+1)^2(j_2+1)^2(j_3+1)^2}
 &\le\frac{c}{(j+1)^2}.
 \label{lemma:tec:convol.2}
\end{align}
\end{lemma}

\begin{proof}The proofs of \eqref{lemma:tec:convol.1} and \eqref{lemma:tec:convol.2} are classical, but we redo them for the sake of completeness. 

To prove \eqref{lemma:tec:convol.1}, consider first the case $1\le l\le j/2$. Then
$j-l+1\ge(j+1)/2$, so the left-hand side of \eqref{lemma:tec:convol.1} is bounded in this case by
$4(j+1)^{-2}(l+1)^{-2}$. The case where $ j/2 \le l \le j-1$ is handled similarly by  using symmetry and the change of variable $\tilde{l}= j-l$. Hence
\begin{equation*}
 \sum_{l=1}^{j-1}\frac{1}{(l+1)^2(j-l+1)^2}
 \le\frac{8}{(j+1)^2}\sum_{l=1}^{\infty}\frac{1}{(l+1)^2},
\end{equation*}
and the series on the right converges, which concludes the proof.

To prove \eqref{lemma:tec:convol.2}, we first sum over $j_1$ and then over $j_2$ and use twice \eqref{lemma:tec:convol.1} to find 
\begin{align*}
\sum_{\substack{j_1+j_2+j_3=j\\j_1,j_2,j_3\ge1}}
 \frac{1}{(j_1+1)^2(j_2+1)^2(j_3+1)^2}&=\sum_{j_2=1}^j\left( \sum_{j_1=1}^{j-j_2-1}\frac1{(j_1+1)^2}\frac1{(j-j_1-j_2+1)^2}\right)\frac1{(j_2+1)^2} \\
& \le c \sum_{j_2=1}^j \frac1{(j-j_2+1)^2}\frac1{(j_2+1)^2} \\ 
& \le c \frac1{(j+1)^2}.
\end{align*}
 
\end{proof}

\begin{proof}[Proof of Proposition \ref{prop:indu:Gammaj}]
To prove \eqref{est:Gammaj_nl}, we prove a slightly stronger property.  There exists $K \ge 1$, $0<\delta \le 1$, such that under the hypothesis \eqref{ass:Ej}, for all $1 \le l \le J$
\begin{equation}\label{est:Gammaj_nl_refined}
 \sup_{0\le t\le1}\left\|\Gamma_t^lu(t)\right\|_{H^1}
 \le \frac{\delta}{(l+1)^2} (KB)^l l!.
\end{equation}

We will prove \eqref{est:Gammaj_nl_refined} by induction. First, we verify it for $l=1$. From the definition of $\Gamma$ in \eqref{Gamma:form.1} and then \eqref{interp:Hk}, it follows that, for all $0 \le t \le 1$, 
\begin{equation*}
\left\|\Gamma_t^1u(t)\right\|_{H^1} \le c \left(\|xu(t)\|_{H^1}+\|\partial_xu(t)\|_{H^1} \right)\le c \|u(t)\|_{\mathcal{H}^2} .
\end{equation*}
From the well-posedness of \eqref{NLS} in $\mathcal{H}^2$ in Proposition \ref{prop:LWP:NLS}, we conclude that
\begin{equation*}
\left\|\Gamma_t^1u(t)\right\|_{H^1} \le  c \|u_0\|_{\mathcal{H}^2} \le \delta KB, 
\end{equation*}
under the condition  \begin{equation} \label{cond:K1}
\delta K \ge  c \|u_0\|_{\mathcal{H}^2} 
\end{equation} 
and recalling $B \ge 1$.
This verifies that \eqref{est:Gammaj_nl_refined} holds in the case $l=1$.

Now, let $2 \le j \le J-1$ be fixed. Assume that \eqref{est:Gammaj_nl_refined} holds for $l=1, \cdots, j-1$. We want to prove that \eqref{est:Gammaj_nl_refined} holds for $l=j$.

Let $0<t\le 1$. By using the Duhamel principle associated to \eqref{equation for Gamma u} with $V=|u|^2$ and \eqref{Gamma:eit}, we have that 
\begin{align} \label{Duhamel:NLS:Gamma}
\Gamma_t^ju(t)=e^{it\partial_x^2}(x^ku_0)
 +i\int_0^te^{i(t-s)\partial_x^2}\Gamma_s^j(|u|^2u)(s)\,ds.
\end{align}
Moreover, formula \eqref{Gamma:form.2} and the Leibniz rule yield 
\begin{align*}
 \Gamma_s^j(|u|^2u)&=e^{\frac{i|x|^2}{4s}}(2is)^{|j|}  \partial_x^j(|\tilde{u}|^2\tilde{u}) \\ &
 =e^{\frac{i|x|^2}{4s}}(2is)^{|j|}\sum_{j_1+j_2+j_3=j}\frac{j!}{(j_1!)(j_2!)(j_3!)}(-1)^c
 (\partial_{x}^{j_1}\tilde{u})(\partial_{x}^{j_2}\tilde{u})\overline{(\partial_{x}^{j_3}\tilde{u})} \\ &
 =\sum_{j_1+j_2+j_3=j}\frac{j!}{(j_1!)(j_2!)(j_3!)}(-1)^{j_3}
 (\Gamma_s^{j_1}u)(\Gamma_s^{j_2}u)(\overline{\Gamma_s^{j_3}u}),
\end{align*}
where we recall the notation $\tilde{u}(x,s)=e^{\frac{i|x|^2}{4s}}u(x,s)$. 
Moreover, we separate the three terms in the above sum for which one of the index $j_i$ is equal to $j$, so that
\begin{equation*}
 \Gamma_s^j(|u|^2u)
 =2|u|^2\Gamma_s^ju+(-1)^ju^2\overline{\Gamma_s^ju}
 +\sum_{\substack{j_1+j_2+j_3=j\\0\le j_1,j_2,j_3\le j-1}}
 \frac{j!}{(j_1!)(j_2!)(j_3!)}(-1)^{j_3}
 (\Gamma_s^{j_1}u)(\Gamma_s^{j_2}u)(\overline{\Gamma_s^{j_3}u}).
\end{equation*}
Thus, we deduce inserting this identity in \eqref{Duhamel:NLS:Gamma}, taking the $H^1$-norm, and then using that $H^1$ is a Banach algebra, that 
\begin{align}
 \|\Gamma_t^ju(t)\|_{H^1}
 &\le \|x^ju_0\|_{H^1}
 +3\int_0^t\|u(s)\|_{H^1}^2\|\Gamma_s^ju(s)\|_{H^1}\,ds\notag\\
&\quad+\sum_{\substack{j_1+j_2+j_3=j\\0\le j_1,j_2,j_3\le j-1}}
 \frac{j!}{(j_1!)(j_2!)(j_3!)}\int_0^t
 \|\Gamma_s^{j_1}u(s)\|_{H^1}
 \|\Gamma_s^{j_2}u(s)\|_{H^1}
 \|\Gamma_s^{j_3}u(s)\|_{H^1}\,ds \notag \\ &=:I+II+III\label{est:Gamma:nls:100}
\end{align}

To estimate the first term on the right-hand side of \eqref{est:Gamma:nls:100}, we use \eqref{est:xjf} and the hypothesis on $E_j$ \eqref{ass:Ej}, 
\begin{align*}
I \le \left(2(j^2+1)E_j(u_0)+2j^2\|u_0\|_{L^2}^2\right)^{\frac12} \le  c\sqrt{j^2+1}B^j j! +cj\|u_0\|_{L^2} .
\end{align*}

From the well-posedness theory in $\mathcal{H}^1$, (see Proposition \ref{prop:LWP:NLS}),
\begin{equation} \label{LWP:NLS:H1}
\|u(s)\|_{H^1} \le c \|u(s)\|_{\mathcal{H}^1} \le c \|u_0\|_{\mathcal{H}^1} ,
\end{equation}
so that we bound the second term on the right-hand side of \eqref{est:Gamma:nls:100} by 
\begin{equation*}
II \le c \|u_0\|_{\mathcal{H}^1}^2\int_0^t\|\Gamma_s^ju(s)\|_{H^1}\,ds .
\end{equation*}

To deal with the last term on the right-hand side  of \eqref{est:Gamma:nls:100}, observe that at least two indexes of $j_1, j_2, j_3$ are larger or equal to $1$. We split the sum depending if the remaining of these indexes is $0$ or larger than $1$. Hence, 
\begin{align*}
III & =3\sum_{l=1}^{j-1}\frac{j!}{l!(j-l)!}\int_0^t \|u(s)\|_{H^1} \|\Gamma_s^{l}u(s)\|_{H^1}
 \|\Gamma_s^{j-l}u(s)\|_{H^1}ds\\  & \quad +\sum_{\substack{j_1+j_2+j_3=j\\1\le j_1,j_2,j_3\le j-1}}
 \frac{j!}{(j_1!)(j_2!)(j_3!)}\int_0^t
 \|\Gamma_s^{j_1}u(s)\|_{H^1}
 \|\Gamma_s^{j_2}u(s)\|_{H^1}
 \|\Gamma_s^{j_3}u(s)\|_{H^1}\,ds.
\end{align*}
Thus, we deduce by using \eqref{LWP:NLS:H1}, the induction hypothesis \eqref{est:Gammaj_nl_refined} for all $1 \le l \le j-1$, and  then  Lemma \ref{lemma:tec:convol}  that
\begin{align*}
III & \le c\|u_0\|_{\mathcal{H}^1}j!\left(\sum_{l=1}^{j-1}\frac1{(l+1)^2(j-l+1)^2} \right)\delta^2 (KB)^j\\  & \quad +cj!\left(\sum_{\substack{j_1+j_2+j_3=j\\1\le j_1,j_2,j_3\le j-1}}
\frac1{(j_1+1)^2(j_2+1)^2(j_3+1)^2}\right)\delta^3 (KB)^j \\ 
& \le c \frac{\left(\|u_0\|_{\mathcal{H}^1}\delta^2+\delta^3\right)}{(j+1)^2} (KB)^j j! 
\end{align*}
By inserting these estimates into \eqref{est:Gamma:nls:100}, we find that 
\begin{align*}
\|\Gamma_t^ju(t)\|_{H^1} \le \sqrt{1+j^2}B^j j! +cj\|u_0\|_{L^2}+ c \frac{\left(\|u_0\|_{\mathcal{H}^1}\delta^2+\delta^3\right)}{(j+1)^2} (KB)^j j!+c \|u_0\|_{\mathcal{H}^1}^2\int_0^t\|\Gamma_s^ju(s)\|_{H^1}\,ds ,
\end{align*}
which implies by using Gronwall's inequality that for all $0\le t \le 1$, 
\begin{align*}
\|\Gamma_t^ju(t)\|_{H^1}  \le e^{c\|u_0\|_{\mathcal{H}^1}^2}\left(\sqrt{1+j^2}B^j j! +j\|u_0\|_{L^2}+ c \frac{\left(\|u_0\|_{\mathcal{H}^1}\delta^2+\delta^3\right)}{(j+1)^2} (KB)^j j!\right) .
\end{align*}
We now fix 
\begin{equation} \label{def:delta}
\delta=\min\left\{\frac1{2ce^{c\|u_0\|_{\mathcal{H}^1}^2}(1+\|u_0\|_{\mathcal{H}^1})},1 \right\} \in (0,1]
\end{equation}
so that 
\begin{equation*}
c e^{c\|u_0\|_{\mathcal{H}^1}^2} \left(\|u_0\|_{\mathcal{H}^1}\delta^2+\delta^3\right) \le \frac{\delta}2 .
\end{equation*}
Finally, we choose $K=K(\|u_0\|_{\mathcal{H}^2}) \ge 1$ satisfying \eqref{cond:K1} and 
\begin{equation*}
c e^{c\|u_0\|_{\mathcal{H}^1}^2} \left(\sqrt{1+j^2} +j\|u_0\|_{L^2}\right) \le \frac{K^j \delta}{2(j+1)^2} .
\end{equation*}
Then, we conclude that 
\begin{equation*}
 \sup_{0\le t\le1}\left\|\Gamma_t^ju(t)\right\|_{H^1}
 \le \frac{\delta}{(j+1)^2} (KB)^j j!,
\end{equation*}
which is \eqref{est:Gammaj_nl_refined} in the case $l=j$. This concludes the proof of \eqref{est:Gammaj_nl_refined} by induction, which implies directly \eqref{est:Gammaj_nl}.

The proof of \eqref{est:Gammaj_nl-} follows similarly by applying the same argument to the function $v(x,t)=u(x,1-t)$ that also solves \eqref{NLS}.
\end{proof}

\subsection{Conclusion of the proof of the unique continuation theorem}\label{sec:conclusion}
In order to prove \eqref{E_pqk:u_0:u_1} we will prove the following proposition. 

\begin{proposition} \label{E_k prop:Ind} Let $0<\alpha, \beta <\frac12$. Let $1<b<\frac1{\beta}-1$ and $1<a<\frac1{\alpha}-1$.
    There exists a constant $K_1 \ge 1$ such that, for all $j \in \mathbb{N},$ $j\ge 1$, we have
    \begin{equation} \label{E_k:ind}
        E_j(u_1) \le K_1^j(j!)^{\frac{2}{b} }
    \end{equation}
    and
        \begin{equation} \label{E_k:ind:u0}
        E_j(u_0) \le K_1^j(j!)^{\frac{2}{a} }
    \end{equation}
\end{proposition}

\begin{proof}
    We will again proceed by induction. For $j=1$, we have by using \eqref{interp:Hk}, and the well-posedness in $\mathcal{H}^2$,
$$E_1(u_1) = \|xu_1\|_{L^2}^2 + \|x\partial_x u_1\|_{L^2}^2 \le c \|u_1\|_{\mathcal{H}^2}^2 \le c\|u_0\|_{\mathcal{H}^2}^2\le K_1 $$ for $K_1>c\|u_0\|_{\mathcal{H}^2}^2.$

Now, let $k\ge 2$ be fixed with $k \ge k_0$, where $k_0$ is given in Proposition \ref{prop:zero-gap}. Assume that \eqref{E_k:ind} holds for $j=1, \cdots, k-1$. We want to prove that \eqref{E_k:ind} holds for $j=k$. 
From \eqref{Almost:monotonicity}, \eqref{comparing:with:tilde1}, and using that $b> 1$, it follows that, for any integer $l \ge 1$,
\begin{align*}
    E_{bl}(u_1) \le c\left(E_{b(l+1)}(\tilde{u}_1) + \|u_1\|_{\mathcal{H}^1}^2\right).
\end{align*}
By Proposition \ref{prop:zero-gap}, with $k=l+1$ it follows that
\begin{align*}
    E_{bl}(u_1)&\le c\left(\|\tilde{u}_1\|_{H^{l+2}}^2 + c^l l^{\frac{2b}{b+1}l}\|u_1\|_{\mathcal{H}^1}^2 \right) \\ &\le c^l \left(\max_{0\le j\le l+1} \|\Gamma^j_1 u_1\|_{H^1}^2 + l^{\frac{2b}{b+1}l   }\|u_1\|_{\mathcal{H}^1}^2\right).
\end{align*}
Fix $b<\tilde{b}<\frac1{\beta}-1$. We choose $l=l(k) = \ceil{\frac{k}{\tilde{b}}}$, \textit{i.e.} $l \in \mathbb N$ and
$\frac{k}{\tilde{b}} \le l < \frac{k}{\tilde{b}}+1$, then $k\le \tilde{b}l$. Using \eqref{Almost:monotonicity} and the above estimate with $b=\tilde{b}$, it follows that
\begin{equation}
    E_{k}(u_1) \le c^l \left(\max_{0\le j\le l+1} \|\Gamma^j_1 u_1\|_{H^1}^2 + l^{\frac{2\tilde{b}}{\tilde{b}+1}l   }\|u_1\|_{\mathcal{H}^1}^2\right).
\end{equation}
Let $k_1=\frac{2\tilde{b}}{\tilde{b}-1}+1$. If $k \ge k_1$, we have $l+1<\frac{k}{\tilde{b}}+2<k$.
The induction hypothesis then yields 
\[ E_j(u_1)\le K_1^j (j!)^{\frac 2{b}} < K_1^j (j!)^2, \quad  \forall  \, 1\le j \le l+1 . \] 
Hence, we can apply Proposition \ref{prop:indu:Gammaj}, with $J=l+1$ and $B=K_1^{\frac12}$.
It follows that there exists $K=K(\|u_0\|_{\mathcal{H}^2}) \ge 1$ such that
\begin{align*} 
     E_{k}(u_1) 
     &\le c^l \left( (KK_1)^{l+1} ((l+1)!)^2 + l^{\frac{2\tilde{b}}{\tilde{b}+1}l   }\|u_1\|_{\mathcal{H}^1}^2 \right) \\ & \le K_1^k \left(c^l K^{l+1} (l!)^2 + l^{\frac{2\tilde{b}}{\tilde{b}+1}l}\|u_1\|_{\mathcal{H}^1}^2\right). 
\end{align*}
Since $l=\frac{k}{\tilde{b}}+O(1)$, it follows by using Stirling's formula that
\begin{align*}
    E_k(u_1) \le K_1^k c^k K^{k+1} (k!)^{\frac2{\tilde{b}}} \le \frac12 K_1^k(k!)^{\frac2{b}} 
\end{align*}
for $k \ge k_2(\|u_0\|_{\mathcal{H}^2})$, where $k_2 \ge \max\{k_0,k_1\}$ is chosen such that
\[  \sup_{k \ge k_1}c^kK^{k+1}e^{2k(\frac1{\tilde{b}}-\frac1{b})\log k} \le \frac12 ,\]
since $b<\tilde{b}$. This proves \eqref{E_k:ind} for $j=k$ and $k\ge k_2.$

For $k<k_2$, we know, since $u_1\in \mathcal{S}(\mathbb R)$ that $E_k(u_1)<\infty$ for all $1 \le k \le k_2$. Then, there exists a constant $K_2\ge 1$ such that 
$$E_k(u_1)\le K_2, \quad \forall \, 1 \le k \le k_2.$$ 
By enlarging the constant $K_1$ if necessary, we conclude that \eqref{E_k:ind} holds for all $k\ge1$.

The proof of \eqref{E_k:ind:u0} follows in a similar way. 
\end{proof} 

Now, we conclude the proof of Theorem \ref{thm:main}.
\begin{proof}[Proof of Theorem \ref{thm:main}]

Let $1<p,q<\infty$ be such that $\frac1p+\frac1q=1$ to be fixed later. By applying \eqref{E_k:ind} for $j=\floor{\frac {qk}{2}}$, and using Stirling's formula, we deduce
\begin{equation}
    E_{\frac{qk}{2}}(u_1) \le C_1^k (k!)^{q\delta_1},
\end{equation} 
for $\delta_1=\frac1b.$
Similarly, \eqref{E_k:ind:u0} with $j=\floor{\frac{pk}{2}}$ yields 
\begin{equation}
    E_{\frac{pk}{2}}(u_0) \le C_1^k (k!)^{p\delta_0}.
\end{equation} 
for $\delta_0=\frac1a.$
Now, we optimize $p$ and $q$ as in the proof of Theorem \ref{theo:UC:potential}. We fix $p$ and $q$ such that $$p\delta_0=q\delta_1 \iff q= \frac{\delta_0+\delta_1}{\delta_1}.$$
Then $p\delta_0<1$ and $q\delta_1<1 \iff \delta_0+\delta_1=\frac1a+\frac{1}{b}<1. $

Remind by hypothesis that $$\frac{\alpha}{1-\alpha}+\frac{\beta}{1-\beta} <1.$$
Now choose $1<a<\frac1{\alpha}-1$ and $1<b<\frac1{\beta}-1$ such that
\[\frac{\alpha}{1-\alpha}+\frac{\beta}{1-\beta}<\frac1a+\frac1b<1, \]
which guarantees then that $p\delta_0<1$ and $q\delta_1<1$.

Therefore, we conclude that
$$\int_{\mathbb R} e^{c|x|^q}|u_1(x)|^2dx =\sum_{k=0}^{\infty} \frac{c^k}{k!} \int|x|^{qk}|u_1(x)|^2 dx \le \sum_{k=0}^{\infty} \frac{c^k}{k!} E_{\frac{qk}{2}}(u_1)\le \sum_{k=0}^{\infty} e^{(q\delta_1-1)k\log{k} + O(k)} <+\infty $$
and 
$$\int_{\mathbb R} e^{c|x|^p}|u_0(x)|^2dx =\sum_{k=0}^{\infty} \frac{c^k}{k!} \int|x|^{pk}|u_0(x)|^2 dx \le \sum_{k=0}^{\infty} \frac{c^k}{k!} E_{\frac{pk}{2}}(u_0) \le \sum_{k=0}^{\infty}e^{(p\delta_0-1)k\log{k} + O(k)} <+\infty $$
and by the the Morgan uncertainty principle (Theorem \ref{theo:EKPV2011}) we deduce that $u\equiv0.$
\end{proof}

\bibliographystyle{acm} 
\bibliography{ref241126}

@article{kulikov_uniqueness,
author = {Kulikov, Aleksei and Nazarov, Fedor and Sodin, Mikhail},
year = {2025},
pages = {84-130},
title = {Fourier Uniqueness and Non-Uniqueness Pairs},
volume = {21},
journal = {Journal of Mathematical Physics, Analysis, Geometry},
doi = {10.15407/mag21.01.04},
}

@article{HasegawaTappert1973,
  author  = {Hasegawa, Akira and Tappert, Frederick},
  title   = {Transmission of stationary nonlinear optical pulses in dispersive dielectric fibers. {I}. {A}nomalous dispersion},
  journal = {Applied Physics Letters},
  volume  = {23},
  number  = {3},
  pages   = {142--144},
  year    = {1973}
}

@book{Agrawal2019,
  author    = {Agrawal, Govind P.},
  title     = {Nonlinear Fiber Optics},
  edition   = {6},
  publisher = {Academic Press},
  year      = {2019}
}

@article{DalfovoGiorginiPitaevskiiStringari1999,
  author  = {Dalfovo, Franco and Giorgini, Stefano and Pitaevskii, Lev P. and Stringari, Sandro},
  title   = {Theory of {B}ose--{E}instein condensation in trapped gases},
  journal = {Reviews of Modern Physics},
  volume  = {71},
  number  = {3},
  pages   = {463--512},
  year    = {1999}
}

@book{PitaevskiiStringari2016,
  author    = {Pitaevskii, Lev and Stringari, Sandro},
  title     = {Bose--{E}instein Condensation and Superfluidity},
  series    = {International Series of Monographs on Physics},
  volume    = {164},
  publisher = {Oxford University Press},
  year      = {2016}
}

@article{BakanHedenmalmMontesRodriguezRadchenkoViazovska2021,
  author  = {Bakan, Andrew and Hedenmalm, Haakan and Montes-Rodr{\'\i}guez, Alfonso and Radchenko, Danylo and Viazovska, Maryna},
  title   = {Fourier uniqueness in even dimensions},
  journal = {Proceedings of the National Academy of Sciences of the USA},
  volume  = {118},
  number  = {15},
  pages   = {e2023227118},
  year    = {2021},
  doi     = {10.1073/pnas.2023227118}
}

@article{Stoller2021,
  author  = {Stoller, Martin},
  title   = {Fourier interpolation from spheres},
  journal = {Transactions of the American Mathematical Society},
  volume  = {374},
  number  = {11},
  pages   = {8045--8079},
  year    = {2021},
  doi     = {10.1090/tran/8440}
}

@article{Kulikov2021,
  author  = {Kulikov, Aleksei},
  title   = {Fourier interpolation and time-frequency localization},
  journal = {Journal of Fourier Analysis and Applications},
  volume  = {27},
  pages = {7},
  paper = {7},
  year    = {2021}
}

@article{CohnKumarMillerRadchenkoViazovska2022,
  author  = {Cohn, Henry and Kumar, Abhinav and Miller, Stephen D. and Radchenko, Danylo and Viazovska, Maryna},
  title   = {Universal optimality of the {$E_8$} and {L}eech lattices and interpolation formulas},
  journal = {Annals of Mathematics},
  volume  = {196},
  number  = {3},
  pages   = {983--1082},
  year    = {2022},
  doi     = {10.4007/annals.2022.196.3.3}
}

@article{RadchenkoStoller2022,
  author  = {Radchenko, Danylo and Stoller, Martin},
  title   = {Fourier non-uniqueness sets from totally real number fields},
  journal = {Commentarii Mathematici Helvetici},
  volume  = {97},
  pages   = {513--553},
  year    = {2022},
  doi     = {10.4171/CMH/538}
}

@article{RamosStoller2022,
  author  = {Ramos, Jo{\~a}o P. G. and Stoller, Martin},
  title   = {Perturbed {F}ourier uniqueness and interpolation results in higher dimensions},
  journal = {Journal of Functional Analysis},
  volume={282},
  number = {12},
  pages   = {109448},
  year    = {2022},
  doi     = {10.1016/j.jfa.2022.109448}
}

@article{BondarenkoRadchenkoSeip2023,
  author  = {Bondarenko, Andriy and Radchenko, Danylo and Seip, Kristian},
  title   = {Fourier interpolation with zeros of zeta and {$L$}-functions},
  journal = {Constructive Approximation},  
  volume={57},
 pages = {405--461},
  year    = {2023}
}

@article{RamosSousa2023,
  author  = {Ramos, Jo{\~a}o P. G. and Sousa, Mateus},
  title   = {Perturbed interpolation formulae and applications},
  journal = {Analysis \& PDE},
  volume  = {16},
  number  = {10},
  pages   = {2327--2384},
  year    = {2023},
  doi     = {10.2140/apde.2023.16.2327}
}

@misc{Adve2023,
  author        = {Adve, Anshul},
  title         = {Density criteria for {F}ourier uniqueness phenomena in {$\mathbb{R}^d$}},
  note        = {Preprint, arXiv:2306.07475 (2023)}
}

@misc{RadchenkoRamos2024,
  author        = {Radchenko, Danylo and Ramos, Jo{\~a}o P. G.},
  title         = {Perturbed lattice crosses and {H}eisenberg uniqueness pairs},
  note = {Preprint, arXiv:2410.04557 (2024)},
  archivePrefix = {arXiv}
}

@misc{RamosSousa2024,
  author        = {Ramos, Jo{\~a}o P. G. and Sousa, Mateus},
  title         = {On {P}auli pairs and {F}ourier uniqueness problems},
  note          = {Preprint, arXiv:2410.12065, (2024)},
  eprint        = {2410.12065},
  archivePrefix = {arXiv}
}

@misc{CasseseRamos2025,
  author        = {Cassese, Gabriele and Ramos, Jo{\~a}o P. G.},
  title         = {On the regularity of {F}ourier interpolation formulas},
  note          = {Preprint, arXiv:2503.15733, (2025)},
  eprint        = {2503.15733},
  archivePrefix = {arXiv}
}

@misc{Lysen2025,
  author        = {Lysen, Torgeir Keun},
  title         = {Critical and asymmetric {F}ourier uniqueness pairs},
  note          = {Preprint, arXiv:2509.17600, (2025)},
  eprint        = {2509.17600},
  archivePrefix = {arXiv}
}

@misc{BerghausBondarenkoRadchenkoSeipSun2025,
  author        = {Berghaus, David and Bondarenko, Andriy and Radchenko, Danylo and Seip, Kristian and Sun, Qihang},
  title         = {The basis functions of {F}ourier interpolation},
  note          = {Preprint, arXiv:2512.18677 (2025)},
  eprint        = {2512.18677},
  archivePrefix = {arXiv}
}

@misc{Lysen2026,
  author        = {Lysen, Torgeir Keun},
  title         = {Discrete {P}auli pairs},
  year          = {2026},
  note = {Preprint, arXiv:2605.23079 (2026)},
  eprint        = {2605.23079},
  archivePrefix = {arXiv}
}

@article{Escauriaza2006Uniqueness,
  author  = {Escauriaza, Luis and Kenig, Carlos E. and Ponce, Gustavo and Vega, Luis},
  title   = {On uniqueness properties of solutions of {S}chr{\"o}dinger equations},
  journal = {Communications in Partial Differential Equations},
  volume  = {31},
  number  = {12},
  pages   = {1811--1823},
  year    = {2006},
  doi     = {10.1080/03605300500530446}
}

@article{Escauriaza2008Convexity,
  author  = {Escauriaza, Luis and Kenig, Carlos E. and Ponce, Gustavo and Vega, Luis},
  title   = {Convexity properties of solutions to the free {S}chr{\"o}dinger equation with {G}aussian decay},
  journal = {Mathematical Research Letters},
  volume  = {15},
  number  = {5},
  pages   = {957--971},
  year    = {2008},
  doi     = {10.4310/MRL.2008.v15.n5.a10}
}

@article{Escauriaza2008Hardy,
  author  = {Escauriaza, Luis and Kenig, Carlos E. and Ponce, Gustavo and Vega, Luis},
  title   = {Hardy's uncertainty principle, convexity and {S}chr{\"o}dinger evolutions},
  journal = {Journal of the European Mathematical Society},
  volume  = {10},
  number  = {4},
  pages   = {883--907},
  year    = {2008},
  doi     = {10.4171/JEMS/134}
}

@article{Escauriaza2010Sharp,
  author  = {Escauriaza, Luis and Kenig, Carlos E. and Ponce, Gustavo and Vega, Luis},
  title   = {The sharp {H}ardy uncertainty principle for {S}chr{\"o}dinger evolutions},
  journal = {Duke Mathematical Journal},
  volume  = {155},
  number  = {1},
  pages   = {163--187},
  year    = {2010},
  doi     = {10.1215/00127094-2010-053}
}

@article{Cowling2010Hardy,
  author  = {Cowling, Michael and Escauriaza, Luis and Kenig, Carlos E. and Ponce, Gustavo and Vega, Luis},
  title   = {The {H}ardy uncertainty principle revisited},
  journal = {Indiana University Mathematics Journal},
  volume  = {59},
  number  = {6},
  pages   = {2007--2025},
  year    = {2010}
}

@article{Escauriaza2011UniqueContinuation,
  author  = {Escauriaza, Luis and Kenig, Carlos E. and Ponce, Gustavo and Vega, Luis},
  title   = {Unique continuation for {S}chr{\"o}dinger evolutions, with applications to profiles of concentration and traveling waves},
  journal = {Communications in Mathematical Physics},
  volume  = {305},
  number  = {2},
  pages   = {487--512},
  year    = {2011},
  doi     = {10.1007/s00220-011-1256-3}
}

@book{Boas1954,
  author    = {Boas, Jr., Ralph Philip},
  title     = {Entire Functions},
  publisher = {Academic Press},
  address   = {New York},
  year      = {1954}
}

@book{GelfandShilov1968vol2,
  author    = {Gel'fand, I. M. and Shilov, G. E.},
  title     = {Generalized Functions. {V}ol. 2: Spaces of Fundamental and Generalized Functions},
  publisher = {Academic Press},
  address   = {New York--London},
  year      = {1968}
}

@article{ChungChungKim1996,
  author  = {Chung, Jaeyoung and Chung, Soon-Yeong and Kim, Dohan},
  title   = {Characterizations of the {G}elfand--{S}hilov spaces via {F}ourier transforms},
  journal = {Proceedings of the American Mathematical Society},
  volume  = {124},
  number  = {7},
  pages   = {2101--2108},
  year    = {1996},
  doi     = {10.1090/S0002-9939-96-03291-1}
}

@article{Boyd2014quartic,
  author  = {Boyd, John P.},
  title   = {The {F}ourier transform of the quartic {G}aussian $\exp(-Ax^4)$: hypergeometric functions, power series, steepest descent asymptotics and hyperasymptotics and extensions to $\exp(-Ax^{2n})$},
  journal = {Applied Mathematics and Computation},
  year    = {2014},
  volume = {241},
  pages = {75--87}
}

@article{Escauriaza2012Survey,
  author  = {Escauriaza, Luis and Kenig, Carlos E. and Ponce, Gustavo and Vega, Luis},
  title   = {Uniqueness properties of solutions to {S}chr{\"o}dinger equations},
  journal = {Bulletin of the American Mathematical Society},
  volume  = {49},
  number  = {3},
  pages   = {415--442},
  year    = {2012},
  doi     = {10.1090/S0273-0979-2011-01368-4}
}

@misc{BondarenkoSeip2026,
  author        = {Bondarenko, Andriy and Seip, Kristian},
  title         = {Fourier-invariant functions with dense zero sets},
  note          = {Preprint, arXiv:2608.13468 (2026)},
  eprint        = {2608.13468},
  archivePrefix = {arXiv}
}

@book{Kirkland2020,
  author    = {Kirkland, Earl J.},
  title     = {Advanced Computing in Electron Microscopy},
  edition   = {3},
  publisher = {Springer},
  address   = {Cham},
  year      = {2020}
}

@article{kehle_uniqueness_2022,
  author       = {Kehle, Christoph and Ramos, João P. G.},
  title        = {Uniqueness of Solutions to Nonlinear Schrödinger Equations from their Zeros},
  journal      = {Annals of PDE},
  year         = {2022},
  volume       = {8},
  number       = {2},
  pages        = {21},
  doi          = {10.1007/s40818-022-00138-1},
  url          = {https://doi.org/10.1007/s40818-022-00138-1},
}

@article{hayashi_nonlinear_1988,
	title = {Nonlinear {Schrödinger} {Equations} in {Weighted} {Sobolev} {Spaces}},
	language = {en},
	journal = {Funkcialaj Ekvacioj},
    volume= {31},
    pages = {363-281},
	author = {Hayashi, Nakao and Nakamitsu, Kuniaki and Tsutsumi, Masayoshi},
	year = {1988},

}

@article{hayashi_solutions_1986,
	title = {On {Solutions} of the {Initial} {Value} {Problem} for the {Nonlinear} {Schrödinger} {Equations} in {One} {Space} {Dimension}},
	language = {en},
	journal = {Math. Z.},
	author = {Hayashi, Nakao and Nakamitsu, Kuniaki and Tsutsumi, Masayoshi},
	year = {1986},
    pages = {637--650},
    volume = {192},
}

@article{hayashi_solutions_1987,
	title = {On {Solutions} of the {Inital} {Value} {Problem} for the {Nonlinear} {Schrödinger} {Equations}},
	journal = {Journal of Functional Analysis},
    volume = {71},
    pages = {218-245},
	author = {Hayashi, Nakao and Nakamitsu, Kuniaki and Tsutsumi, Masayoshi},
	year = {1987},
}

@article{ramos_fourier_2022,
	title = {Fourier uniqueness pairs of powers of integers},
	volume = {24},
	issn = {1435-9855, 1435-9863},
	url = {https://ems.press/doi/10.4171/jems/1194},
	doi = {10.4171/jems/1194},
	language = {en},
	number = {12},
	urldate = {2024-09-09},
	journal = {Journal of the European Mathematical Society},
	author = {Ramos, João P. G. and Sousa, Mateus},
	month = apr,
	year = {2022},
	pages = {4327--4351},
}

@article{Escauriaza2011Morgan,
  author  = {Escauriaza, Luis and Kenig, Carlos E. and Ponce, Gustavo and Vega, Luis},
  title   = {Uncertainty principle of Morgan type and Schr{\"o}dinger evolutions},
  journal = {Journal of the London Mathematical Society},
  volume   = {83},
  number   = {1},
  pages    = {187--207},
  year     = {2011},
  doi      = {10.1112/jlms/jdq072}
}

@book{rodino1993linear,
  title={Linear Partial Differential Operators in Gevrey Spaces},
  author={Rodino, L.},
  isbn={9789810208455},
  lccn={96142627},
  series={G - Reference,Information and Interdisciplinary Subjects Series},
  url={https://books.google.no/books?id=02bcuV7ObvYC},
  year={1993},
  publisher={World Scientific}
}

@book{Hormander2003AnalysisI,
  author    = {H{\"o}rmander, Lars},
  title     = {The Analysis of Linear Partial Differential Operators I:
               Distribution Theory and Fourier Analysis},
  series    = {Classics in Mathematics},
  publisher = {Springer-Verlag},
  address   = {Berlin},
  year      = {2003},
  isbn      = {978-3-540-00662-1},
  doi       = {10.1007/978-3-642-61497-2}
}

@article{RadchenkoViazovska2019,
author = {Danylo Radchenko and Maryna Viazovska},
title = {Fourier interpolation on the real line},
journal = {Publications Mathématiques de l'IHÉS},
volume = {129},
pages = {51--81},
year = {2019},
doi = {10.1007/s10240-018-0101-z},
publisher = {Springer Berlin Heidelberg}
}

\end{document}